\documentclass{tran-l}
\usepackage{amsmath}
\usepackage{amsfonts}
\usepackage{amssymb}
\usepackage{epsfig}
\usepackage{psfrag}

\newcommand{\begriff}[1]{\textbf{#1}}

\renewcommand{\d}{\partial}
\newcommand{\nach}{\rightarrow}
\newcommand{\T}{\mathcal T}
\newcommand{\C}{\mathcal C}

\newcommand{\z}{\mathcal Z}

\newcommand{\co}{\colon\thinspace}    %  Colon with correct spacing for maps.

\newtheorem{theorem}{Theorem}[section]
\newtheorem{lemma}[theorem]{Lemma}
                                  
\theoremstyle{definition}
\newtheorem{definition}[theorem]{Definition}

\newtheorem{construction}[theorem]{Construction}

\theoremstyle{remark}

\numberwithin{equation}{section}

\begin{document}

\title{How to make a Triangulation of $S^3$ polytopal}  

\author{Simon A. King}
\address{Department of Mathematics,
  Darmstadt University of Technology,
  Schlossgartenstr.~7,
  64289 Darmstadt,
  Germany}
\email{king@mathematik.tu-darmstadt.de}

\begin{abstract}
  We introduce a numerical isomorphism invariant $p(\T)$ for
  any triangulation $\T$ of $S^3$. Although its definition is purely
  topological (inspired by the bridge number of knots), $p(\T)$
  reflects the geometric properties of $\T$.
  Specifically, if $\T$ is polytopal or shellable then $p(\T)$ is
  ``small'' in the sense that we obtain a linear upper bound for
  $p(\T)$ in the number $n=n(\T)$ of tetrahedra of $\T$.
  Conversely, if $p(\T)$ is ``small'' then $\T$ is ``almost''
  polytopal, since we show how to transform $\T$ into a polytopal
  triangulation by $O((p(\T))^2)$ local subdivisions.
  The minimal number of local subdivisions needed to transform $\T$
  into a polytopal triangulation is at least $\frac{p(\T)}{3n}-n-2$.
  
  Using our previous results [The size of triangulations supporting a
  given link.  {\em Geometry \& Topology} \textbf{5} (2001),
  369--398], we obtain a general upper bound for $p(\T)$ exponential
  in $n^2$.  We prove here by explicit constructions that there is no
  general subexponential upper bound for $p(\T)$ in $n$. Thus, we
  obtain triangulations that are ``very far'' from being polytopal.
  
  Our results yield a recognition algorithm for $S^3$ that is
  conceptually simpler, though somewhat slower, as the famous
  Rubinstein--Thompson algorithm.
\end{abstract}

% 57Q15 Triangulating manifolds
% 52B11 $n$-dimensional polytopes
% 52B22 Shellability
% 52B55 Computational aspects related to convexity 
% 57M15 Relations [of low dim. Topology] with graph theory 
% 05C10 Topological graph theory, imbedding 
% 57M25 Knots and links in $S^3$ 
% 57M27 Invariants of knots and 3-manifolds 

\subjclass[2000]{Primary 52B11, 57M25; Secondary 57M15, 05C10, 52B22}

\keywords{Convex polytope, dual graph, spatial graph, polytopality,
bridge number, recognition of the $3$--sphere}

\maketitle

%**********************************************************************

\section{Introduction}
\label{sec:theorems}

\subsection{Method and results}

A cellular decomposition of the $d$--dimensional sphere $S^d$ is
\begriff{polytopal} if it is isomorphic to the boundary complex of a
convex $(d+1)$--polytope.  The study of polytopal cellular
decompositions has a long history and is still an important branch of
research.
By a theorem of Steinitz~\cite{steinitz} of 1922, all triangulations
of $S^2$ are polytopal. 
However, ``most'' triangulations of higher-dimensional spheres are
not polytopal: Kalai~\cite{kalai} has shown that the number of
triangulations of $S^d$ grows faster in the number of vertices than
the number of {polytopal} triangulations of $S^d$, for $d\ge 4$.  
Recently, Pfeifle and Ziegler~\cite{pfeifleziegler} have obtained a
similar result for triangulations of $S^3$.  Thus, for a thorough
understanding of triangulations of spheres, one has to go beyond the
geometric setting of convex polytopes. 

\begin{figure}
  \begin{center}
    {\psfrag{e}{$e$}\epsfig{file=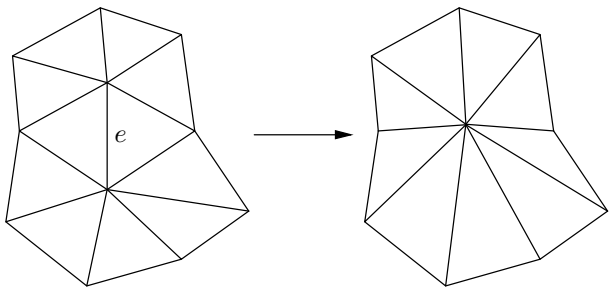}}
    \caption{Contraction along an edge, in dimension $2$}
    \label{fig:contraction}
  \end{center}
\end{figure}

This paper is concerned with triangulations $\T$ of the
$3$--dimensional sphere $S^3$. We introduce an invariant $d(\T)$ of
$\T$ that measures to what extent $\T$ fails to be polytopal, and an
invariant $p(\T)$ that measures ``knottedness'' of $\T$.
We establish a close relationship of $d(\T)$ and $p(\T)$, thus a
relationship of geometric and topological properties of $\T$.
The definition of $d(\T)$ is based on the following local transformations; compare~Figure~\ref{fig:contraction}.
\begin{definition}
  Let $M$ be a closed PL-manifold with 
  PL-triangulations $\T_1$ and $\T_2$, and let  $e$ be an edge of
  $\T_1$ with $\d e=\{a,b\}$. Suppose that $\T_2$ is obtained from
  $\T_1$ by removing the open star of $e$ and identifying the join $a*\sigma$
  with $b*\sigma$ for any simplex $\sigma$ in the link of $e$.
  Then $\T_2$ is obtained from $\T_1$ by the \begriff{contraction}
  along $e$, and $\T_1$ is obtained from $\T_2$ by an
  \begriff{expansion} along $e$. 
\end{definition}
In general, there are edges of $\T_1$ along which the contraction is 
impossible. This is the case, e.g., if an edge $e$ of $\T_1$ is part
of an edge path of length $3$ that does not bound a $2$--simplex of
$\T_1$. Indeed then $\T_2$ has multiple edges and is not a simplicial
complex.
Obviously any PL--triangulation admits only a finite number of
contractions. 

Stellar subdivisions of simplices are examples of expansions. 
Since an expansion increases the number of vertices by
one and the number of simplicial complexes with a given number of
vertices is finite, it is easy to see that also the number of possible
expansions is finite. 

By results of Alexander~\cite{alexander} and Moise~\cite{moise}, any
two triangulations of a $3$--manifold are related by a finite sequence
of stellar subdivisions and their inverses.  From this, one can
conclude that any triangulation $\T$ of $S^3$ can be transformed into
a polytopal triangulation by a finite sequence of expansions.  We
define $d(\T)$ as the length of a shortest sequence of expansions
relating $\T$ with a polytopal triangulation. It is a notion of the
geometric complexity of $\T$.

Our basic idea is to study $\T$ via its embedded dual graph
$\C^1\subset S^3$, i.e., the $1$--skeleton of the dual
cellular decomposition $\C=\T^*$ of $\T$.
The embedding of $\C^1$ is determined by $\T$ up to ambient isotopy. 
Thus, any ambient isotopy invariant of $\C^1\subset S^3$ is an 
(oriented) isomorphism invariant of $\T$. This allows to use techniques
from knot theory in discrete geometry.
Note that if $\T$ is polytopal then its isomorphism class is already
determined by its \emph{abstract} unembedded dual graph;
see~\cite{ziegler}.  However, there are non-polytopal triangulations
with isomorphic abstract dual graphs whose embedded dual graphs are
not ambiently isotopic to each other.

As the first example of an ambient isotopy invariant of
the dual graph of $\T$, we introduce the \emph{polytopality} $p(\T)$ of $\T$.
We outline its definition.
We denote by $\C^i$ the $i$-skeleton of $\C$, for $i=0,\dots,3$.  Let
$H:S^2\times [0,1]\to S^3$ be an embedding in general position to
$\C$. A parameter $\xi_0\in [0,1]$ is a \emph{critical parameter} of $H$
with respect to $\C^i$, if the surface $H(S^2\times \xi_0)$ is not
in general position to $\C^i$. The number of critical parameters of $H$ with
respect to $\C^i$ is denoted by $c(H,\C^i)$.
The \emph{polytopality} of $\T$ is defined as $p(\T)=\min_H
c(H,\C^1)$, where the minimum is taken over all embeddings
$H:S^2\times [0,1]\to S^3$ in general position to $\C$ with
$\C^2\subset H(S^2\times [0,1])$. 

This definition generalizes the bridge number of knots (introduced by
Schubert~\cite{schubert1} in 1954) to spatial graphs, up to a factor
$2$. The bridge number of knots formed by edges of $\T$ has been
studied earlier, for instance, by Lickorish, Armentrout, Ehrenborg and
Hachimori: One obtains linear upper bounds in $n$ for the bridge
number of such knots, provided $\T$ is shellable~\cite{lickorish}, its
dual is shellable~\cite{armentrout} or $\T$ is vertex
decomposable~\cite{ehrenborghachimori}.  The following theorem is in
the same spirit.  For the notions of diagrams and of shellable cell
complexes, see~\cite{ziegler}.

\begin{theorem}\label{thm:poly3}
  Let $\T$ be a triangulation of $S^3$ with $n$ tetrahedra. 
  \begin{enumerate}
  \item If $\T$ is polytopal, then $p(\T)=n$.
  \item If $\T$ has a diagram, then $p(\T)\le 3n$.
  \item If $\T$ or its dual is shellable, then $p(\T)\le 7n$.
  \end{enumerate}
\end{theorem}

It is natural to ask what happens if one drops the geometric
assumptions on $\T$. It turns out that the estimates for $p(\T)$, and
similarly for the bridge number of links formed by edges of $\T$,
change dramatically, by the following result.
\begin{theorem}\label{thm:poly2}
  For any $m\in \mathbb N$ there is a triangulation $\T_m$ of $S^3$
  with at most $856 m + 534$ tetrahedra and $p(\T_m) > 2^{m-1}$.
\end{theorem}
Hence there is no general subexponential upper
bound for the polytopality of a triangulation, in terms of the number
of tetrahedra. We obtain the following general bound from our
results in~\cite{king1} that are based on a study of the
Rubinstein--Thompson algorithm.
\begin{theorem}\label{thm:poly1}
  If $\T$ is a triangulation of $S^3$ with $n$ tetrahedra, then
  $$n\le p(\T)<2^{200 n^2}.$$
\end{theorem}

The definition of $p(\T)$ does not rely on geometry. Nevertheless, it turns
out that $p(\T)$ is in close relationship with geometrical properties
of $\T$, as we obtain both lower and upper bounds for $d(\T)$
in terms of $p(\T)$. 
We start with a lower bound.
\begin{theorem}\label{thm:basic3b}
  If $\T$ is a triangulation of $S^3$ with $n$ tetrahedra, then
  $$d(\T) > \frac{p(\T)}{2n+1} - n - \frac 53.$$ 
\end{theorem}
This together with Theorem~\ref{thm:poly2} implies that
there are triangulations of $S^3$ that are ``very far'' from
being polytopal. 
So in a certain sense it reflects Pfeifle's and Ziegler's
result~\cite{pfeifleziegler} that ``most'' triangulations of $S^3$ are
not polytopal.

To obtain upper bounds for $d(\T)$ in $p(\T)$, we consider the class
of \begriff{edge contrac\-tible} triangulations. A triangulation of
$S^d$ is edge contractible if one can transform it into the boundary
complex of a $(d+1)$--simplex by successive contractions along edges.
It is well known and easy to show that any edge contractible
triangulation is polytopal.
All triangulations of $S^2$ are edge contractible by a theorem of
Wagner~\cite{wagner}.  In~\cite{klee}, Section~6.3, one finds an
example of a polytopal triangulation of $S^3$ that is not edge
contractible.

\begin{theorem}\label{thm:basic3}
  From any triangulation $\T$ of $S^3$, one can
  obtain an edge contractible triangulation of $S^3$ by a sequence of at
  most $512 (p(\T))^2 + 869 p(\T) + 376$ successive expansions.  In
  particular, if $n$ is the number of tetrahedra,
  $$d(\T) \le  512 (p(\T))^2 + 869 p(\T) + 376 < 2^{401 n^2}.$$
\end{theorem}

The bounds in Theorem~\ref{thm:basic3} are certainly not sharp. 
Theorem~\ref{thm:poly3} implies that under additional
geometric or combinatorial assumptions on $\T$ one can replace the
estimate of Theorem~\ref{thm:basic3} by a quadratic bound in $n$.
However, a general upper bound for $d(\T)$ must be at least
exponential in $n$, by Theorems~\ref{thm:poly2}
and~\ref{thm:basic3b}. 
Recently we applied Theorem~\ref{thm:basic3} to prove an upper
bound for the \emph{crossing number} of links embedded in the $1$--skeleton
of $\T$; see~\cite{kingcross}.

According to Theorems~\ref{thm:poly3} and~\ref{thm:basic3b}, if $\T$
is geometrically ``simple'' then $p(\T)$ is small.
According to Theorem~\ref{thm:basic3}, if $p(\T)$ is small then
$\T$ is not far from being polytopal. 
Thus $p(\T)$ turns out to be a measure
for the geometric complexity of $\T$.

As an intermediate step in the proof of Theorem~\ref{thm:basic3}, we
establish the following linear upper bound for the ``distance'' of
two triangulations of $S^3$.
\begin{theorem}\label{thm:basic2}
  Any two triangulations $\T,\widetilde \T$ of $S^3$ 
  are related by a sequence of at most $325 \left(p(\T) + 
  p(\widetilde \T)\right) + 508$ expansions and  contractions.
\end{theorem}
This together with Theorem~\ref{thm:poly1} implies that any two
triangulations of $S^3$ with $\le n$ tetrahedra can be related by a
sequence of less than $2^{201 n^2}$ contractions and expansions.
Mijatovi\'c~\cite{mijatovic} obtained a similar estimate concerning
bistellar moves rather than contractions and expansions, using a
modified form of our results in~\cite{king1}.
Under additional geometric or combinatorial assumptions on $\T$,
Theorem~\ref{thm:basic2} together with
Theorem~\ref{thm:poly3} yields a linear bound in $n$.
Recently we have extended our method to the projective space
$P^3_{\mathbb R}$, proving that any two triangulations of
$P^3_{\mathbb R}$ with at most $n$ tetrahedra can be related by a
sequence of less than $2^{27000\,n^2}$ edge contractions and
expansions~\cite{kingproj}.

Theorem~\ref{thm:basic2} yields the following simple try-and-check
algorithm for the recognition of $S^3$. Let $\T$ be a triangulation of
a closed $3$--manifold $N$ with $n$ tetrahedra. Try all sequences of $<
2^{201 n^2}$ contractions and expansions starting from $\T$, which are
finite in number. If one of them transforms $\T$ into the boundary
complex of a $4$--simplex then $N$ is homeomorphic to $S^3$. Otherwise,
$N$ is not homeomorphic to $S^3$ by Theorem~\ref{thm:basic2}.
This recognition algorithm is certainly slower, but simpler,
as the Rubinstein--Thompson recognition algorithm.

%------------------------------

\subsection{Outline and organization of the proofs}

We outline the proof of Theorem~\ref{thm:poly3}.
Let $\T$ be a triangulation of $S^3$, and let 
$\C$ be its dual cellular decomposition.
If $\T$ is polytopal (resp. $\T$ has a diagram), then $\C$ 
(resp. its barycentric subdivision $\C'$) has a diagram. 
A sweep-out of $\mathbb R^3$ by Euclidean planes in general
position to the diagram of $\C$ (resp. $\C'$) has only critical points
in the vertices of the diagram. Hence $p(\T)$ is bounded by the number
of vertices (resp. the number of vertices plus the number of edges) of
$\C$. 
If $\T$ or $\C$ is shellable, then $\C'$ is shellable. Using a
shelling order of the tetrahedra of $\C'$, we construct an embedding
$H\co S^2\times I\to S^3$ that
has at most one critical point in each open simplex of $\C'$. This
yields a bound for $p(\T)$.

Theorem~\ref{thm:poly2} is based
on a result of Hass, Snoeyink and Thurston~\cite{hass}.
We construct for any $m\in \mathbb N$ a simple cellular decomposition
$\z_m$ of the solid torus with a linear bound in $m$ for the number of
vertices, so that any meridional disc for the torus intersects
$\z_m^1$ in $\ge 2^{m-1}$ points. We glue two copies of $\z_m$
together and obtain a simple cellular decomposition $\C_m$ of $S^3$
that is dual to a triangulation $\T_m$ of $S^3$. 
It follows $p(\T_m)> 2^{m-1}$, since for any sweep-out there is a
level sphere that contains a meridional disc in one of the solid tori.

We outline the proof of Theorem~\ref{thm:basic3b}.  Let $\T$ be a
triangulation of $S^3$ with $n$ tetrahedra.  Let $\tilde \T$ be a
polytopal triangulation that is obtained from $\T$ by a sequence of
$d(\T)$ expansions.  It has at most $n+ d(\T)$ vertices.  Since
$\tilde \T$ is polytopal, there is an embedding $H\co S^2\times I \to
S^3$ with $\tilde\T^2\subset H(S^2\times I)$ that has only critical
points in the vertices of $\tilde\T$, i.e., $c(H,\tilde\T^1)\le
n+d(\T)$.
We have $c(H,\T^1)\le c(H,\tilde\T^1)$, since $\T^1\subset \tilde\T^1$.
We finish the proof of Theorem~\ref{thm:basic3b} by giving a lower bound
for $c(H,\T^1)$ in terms of $p(\T)$ and $n$.
This lower bound also yields Theorem~\ref{thm:poly1}  from our results 
in~\cite{king1}.

The proofs of Theorems~\ref{thm:basic3} and~\ref{thm:basic2} are based on
the interplay of cellular structures on $S^3$ and isotopies of surfaces.
Let $\T$ be a triangulation of $S^3$, and let $\C$ be its dual
cellular decomposition. 
Let $H\co S^2\times I\to S^3$ be an embedding in general position
with respect to $\C$ such that $\C^2\subset H(S^2\times I)$ and
$c(H,\C^1)=p(\T)$. 
In the first step of the proof of Theorem~\ref{thm:basic2}, we change
$H$ by canceling pairs of critical parameters so that
$c(H,\C^2)$ is bounded in terms of $p(\T)$.
We associate to any non-critical parameter $\xi\in I$ a $2$--dimensional
polyhedron 
$P_\xi \subset S^3$ that consists of $H(S^2\times \xi)$ and parts of 
$\C^2$.  
It changes by insertions and deletions of $2$--strata,
when $\xi$ passes a critical 
parameter of $H$ with respect to $\C^2$. 
We can bound the number of vertices of the inserted or deleted $2$--strata.

The polyhedron $P_\xi$ is the $2$--skeleton of a cellular decomposition
of $S^3$, whose barycentric subdivision is a triangulation $\T_\xi$.
The deletion of a $2$--stratum $c$ in $P_\xi$ (resp.\ its insertion 
into $P_\xi$) gives rise to a sequence of contractions (resp.\ expansions) 
starting at $\T_\xi$. 
The number of contractions (resp.\ expansions) is  determined by the
number of vertices in $\d c$. 
We obtain a sequence of expansions and 
contractions that relates $\T$ with the
boundary complex of a $4$--simplex, whose length is bounded in terms of
$p(\T)$.  
This yields Theorem~\ref{thm:basic2}.

The proof idea for Theorem~\ref{thm:basic3} is to attach $2$--strata
to $P_\xi$ as in the proof of Theorem~\ref{thm:basic2}, when $\xi$
passes a critical parameter of $H$, and to 
postpone the deletions of $2$--strata until all insertions are done. 
This gives rise to a sequence of
expansions followed by a sequence of contractions that relate $\T$ with the
boundary complex of a $4$--simplex. Here the boundary of the inserted
$2$--strata may contain additional vertices, at the points of
intersections with the boundary of other inserted (and not yet deleted)
$2$--strata. This yields the quadratic bound in
Theorem~\ref{thm:basic3}. The bound in terms of $n$ is a consequence of
Theorem~\ref{thm:poly1}.  

%----------------------------------------------

The paper is organized as follows. 
In Subsection~\ref{sec:defpoly}, we define the polytopality and prove
Theorems~\ref{thm:poly3}, \ref{thm:poly1} and~\ref{thm:basic3b}. 
Subsection~\ref{sec:expo} is devoted to the proof of
Theorem~\ref{thm:poly2}. 
Theorems~\ref{thm:basic3} and~\ref{thm:basic2} are proven in
Section~\ref{sec:hauptbeweise}. 
In Subsection~\ref{sec:minmax}, 
we change any embedding $H\co
S^2\times I\to S^3$ into an embedding $\tilde H\co S^2\times I\to S^3$
with an upper bound for $c(\tilde H,\C^2)$ in
terms of $c(H,\C^1)= c(\tilde H,\C^1)$. 
The proofs of Theorems~\ref{thm:basic2} and~\ref{thm:basic3} are
finished in Subsections~\ref{sec:proof1} and~\ref{sec:proof2}, 
respectively.

%***********************************************************************

\section{Polytopality}
\label{sec:poly}

\subsection{Definition and upper bounds}
\label{sec:defpoly}

After exposing some technical notions, we introduce in this subsection a
numerical invariant for a triangulation $\T$ of $S^3$, called
\emph{polytopality}.  
It is an ambient isotopy invariant of the embedded dual graph of $\T$ and
is inspired by the bridge number 
of links due to Schubert~\cite{schubert1}.
We prove a general upper bound for $p(\T)$ in terms of the number of
tetrahedra, based on our results in~\cite{king1}. Under different
hypotheses on $\T$, we obtain stronger bounds. We prove a lower bound
for $d(\T)$ in terms of $p(\T)$.

Let $M$ be a closed $3$--manifold. 
For a cellular decomposition $\z$ of $M$ and for $i=0,\dots, 3$, let
$\z^i$ denote the $i$--skeleton of $\z$.
In this paper we only consider cellular decompositions
such that if $\varphi\co B^{k}\to M$ is the attaching map
of a $k$--cell, then for any open cell $c$ the restriction of
$\varphi$ to a connected component of $\varphi^{-1}(c)$ is a
homeomorphism onto $c$.  
We do not assume that the closure of an open $k$--cell is a compact
$k$--ball. 
Recall that a \begriff{triangulation} of $M$ is
a cellular decomposition of $M$ that forms a simplicial complex.
We denote by $\#(X)$ the number of connected components of a
topological space $X$. The notation $X\subset M$ stands for a tame
embedding of $X$ into $M$, and $U(X)$ denotes an open regular
neighborhood of $X$ in $M$. We denote the unit interval by $I=[0,1]$.

Let $\z\subset M$ be a $2$--dimensional cell complex (later, it will be the
$1$-- or $2$--skeleton of a triangulation of $M$ or of its dual).  An
\begriff{isotopy mod $\z$} is an ambient isotopy that preserves each
open cell of $\z$ as a set.
Let $S$ be a closed surface. 
Let $H\co S\times I\to M$ be an embedding.  For $\xi\in I$, set
$H_\xi=H(S\times \xi)$. A number $\xi\in I$ is a \begriff{critical
  parameter} of $H$ with respect to $\z$, and a point $p\in H_\xi$ is
a \begriff{critical point} of $H$ with respect
to $\z$, if $p$ is a vertex of $\z$, a point of tangency of $H_\xi$ with
$\z^1$, or a point of tangency of $H_\xi$ with $\z^2$.
An embedding $H\co S\times I\to M$ is \begriff{$\z^1$--Morse} if it
has finitely many critical parameters with respect to $\z$, to each
critical parameter of $H$ with respect to $\z$ belongs exactly one
critical point, and one connected component of $U(p_0)\setminus
H_{\xi_0}$ is disjoint from $\z^1$, for each critical point $p_0\in
\z^1\setminus \z^0$.
  If $H$ is a $\z^1$--Morse embedding then $c(H,\mathcal \z^i)$
  denotes the number of critical points of $H$ in $\z^i$, for $i=1,2$.

\begin{definition}
  Let $\T$ be a triangulation of $S^3$, and let $\C$ be its dual
  cellular 
  decomposition. The \begriff{polytopality} of 
  $\T$ is the number $$p(\T) = \min_H \,c(H,\C^1),$$ where the minimum is
  taken over all $\C^1$--Morse embeddings $H\co S^2\times I\to S^3$ with
  $\C^2\subset H(S^2\times I)$.
\end{definition}

The following lemma allows to reduce Theorem~\ref{thm:poly1} to 
our results in~\cite{king1}, and is also useful to prove 
Theorem~\ref{thm:basic3b}.
\begin{lemma}\label{lem:T2C2}
  Let $\T$ be a triangulation of $S^3$ with $n$ tetrahedra, and let
  $H\co S^2\times I\to S^3$ be a $\T^1$--Morse embedding with $\T^2\subset
  H(S^2\times I)$. Then $$p(\T) < (2n+1)\cdot
  c(H,\T^1) + 5n.$$
\end{lemma}
\begin{proof}
  Let $\T'$ be the barycentric subdivision
  of $\T$. All regular 
  neighborhoods occuring in this proof are to understand with respect
  to $\T$.
  For any simplex $\sigma$ of $\T$, choose a vertex $v_\sigma$ of
  $\sigma$. 
  Let $p_\sigma\in (\T')^0$ be the barycenter of $\sigma$. 
  By ambient isotopy of $(\T')^1$ with support in $U(\sigma)$, we can
  assume that $p_\sigma\in  U(v_\sigma)$. 

  Let $\tau$ be a  boundary simplex of $\sigma$, and let $e$ be the
  edge of $\T'$ with endpoints $p_\sigma$, $p_\tau$. 
  If $v_\sigma=v_\tau$ then
  we can assume by ambient isotopy of $(\T')^1$ with support in
  $U(\sigma)$ that $e\subset U(v_\sigma)$ and
  $H$ has no critical points in the interior of $e$.
  If $v_\sigma\not=v_\tau$ then let $f$ be the edge of $\T$ with
  endpoints $v_\sigma,v_\tau\in \d \sigma$. We can assume by ambient
  isotopy of $(\T')^1$ with support in 
  $U(\sigma)$ that $e\subset U(f)$ and that the critical points of $H$ in the
  interior of $e$ are in bijective correspondence to those in the
  interior of $f$.
  The latter case occurs at most $2n+1$ times
  for any edge of $\T$, since the star of $e$ in $\T$ contains $\le
  2n+1$ simplices. 
  Thus $H$ has $\le (2n+1)\cdot (c(H,\T^1)-
  \#(\T^0))$ critical points in $(\T')^1\setminus (\T')^0$.
  The $1$--skeleton of the dual cellular decomposition of $\T$ is
  contained in $(\T')^1$. Thus 
  $$p(\T) \le c(H,(\T')^1) \le  (2n+1)\cdot \left(c(H,\T^1)-
  \#(\T^0)\right) + \#\left((\T')^0\setminus \T^0\right).$$ 
  Since $\T$ has $n$ tetrahedra, $2n$ $2$--simplices and at most $2n$
  edges, $(\T')^0\setminus \T^0$ comprises $\le 5 n $ vertices.
\end{proof}

\medskip
\noindent
\emph{Proof of Theorem~\ref{thm:poly1}.}
  Let $\T$ be a triangulation
  of $S^3$ with $n$ tetrahedra and let $\C$ be its dual cellular
  decomposition. 
  Since the vertices of $\C$ are critical points, it
  follows $p(\T)\ge n$.  
  By our results in~\cite{king1}, that we have slightly improved 
  in Chapters~3 and 4 of~\cite{king2}, there is a $\T^1$--Morse 
  embedding $H\co S^2\times I\to S^3$ with $\C^2\subset H(S^2\times I)$ 
  and $c(H,\T^1) < 2^{190 n^2}$. By the preceding Lemma~\ref{lem:T2C2}, 
  we have 
  $$p(\T) < (2n + 1)\cdot c(H,\T^1) + 5 n < 2^{200 n^2}.\qed$$

\medskip
\noindent
\emph{Proof of Theorem~\ref{thm:basic3b}.}
  Let $\T$ be a triangulation of $S^3$ with $n$ tetrahedra, and 
  let $\tilde\T$ be a polytopal
  triangulation with $\tilde v$ vertices that is obtained from $\T$ 
  by $d(\T)$ expansions.
  Since $\tilde\T$ has a diagram, there is a $\tilde\T^1$--Morse embedding
  $H\co S^2\times I\to S^3$ such that $\tilde \T^1\subset H(S^2\times I)$ and
  $c(H,\tilde\T^1)$ equals the number of vertices of $\tilde\T$.
  Since $\T^1\subset \tilde\T^1$ and by Lemma~\ref{lem:T2C2}, we have
  $$  \tilde v = c(H,\tilde\T^1) \ge c(H,\T^1) >
             \frac{p(\T) - 5n}{2n+1}.$$
  The number $v$ of vertices of $\T$ is bounded from above by $n$. 
  Since $d(\T)= \tilde v - v$, we obtain
  $$  d(\T) > \frac{p(\T)}{2n+1} - n - \frac 53.\qed$$
\medskip

The rest of this section is devoted to the proof of Theorem~\ref{thm:poly3}.
The three separate claims are proved in the following three lemmas.  
Recall that a \begriff{$d$--diagram} is  a decomposition 
of a convex $d$--polytope into 
convex polytopes; see~\cite{ziegler} for details. 
A cellular decomposition of $S^{d+1}$ has a diagram, if by removing one
of its top-dimensional cells it becomes isomorphic to a $d$--diagram.
It is well known that any polytopal cellular 
decomposition of $S^d$ has a so-called
\emph{Schlegel diagram}
(named after Schlegel~\cite{schlegel}).  

\begin{lemma}\label{lem:poly-pt}
  Let $\T$ be a triangulation $\T$ of $S^3$ with $n$ tetrahedra. If $\T$ is polytopal then $p(\T)=n$.
\end{lemma}
\begin{proof}
Let $\C$ be the dual cellular decomposition of $\T$. Since $\T$ is
polytopal, $\C$ is polytopal as well. Thus $\C$ has a
Schlegel diagram $\mathcal D\subset \mathbb R^3=S^3\setminus
\{\infty\}$. 
We choose coordinates $(x,y,z)$ for $\mathbb R^3$ so that no
edge of $\mathcal D$ is parallel to the $xy$--plane. 
Then a sweep-out of $\mathbb R^3$ by
planes parallel to the $xy$--plane
gives rise to a $\C^1$--Morse embedding having only critical points in
the $n$ vertices of $\mathcal C$, i.e., $p(\T)=n$.  
\end{proof}
\begin{lemma}
  Let $\T$ be a triangulation $\T$ of $S^3$ with $n$ tetrahedra. If
  $\T$ has a diagram then $p(\T) \le 3 n$.
\end{lemma}
\begin{proof}
  Let $\T'$ be the barycentric subdivision of $\T$. Let $\Gamma\subset
  (\T')^1$ be 
  the $1$--skeleton of the dual cellular decomposition of $\T$. 
  Since $\T$ has a diagram, $\T'$ also has a diagram. A sweep-out by
  Euclidean planes yields
  a $(\T')^1$--Morse embedding $H\co S^2\times I\to S^3$ with $(\T')^0$ as
  set of critical points. The critical points of $H$ in $\Gamma$ are
  thus the $n$ vertices of $\Gamma$ and at most one point in each of
  the $2n$ open edges of $\Gamma$. This yields the lemma. 
\end{proof}
\begin{lemma}
  Let $\T$ be a triangulation  of $S^3$ with $n$ tetrahedra. If $\T$
  or its dual is shellable, then  $p(\T) \le 7 n$.
\end{lemma}
\begin{proof}
  Let $\T'$ be the barycentric subdivision of $\T$. 
  Since $\T$ or its dual is shellable, also $\T'$ is shellable. 
  Thus there is a shelling order $t_1,\dots, t_{24n}$ on the open
  tetrahedra of $\T'$ such that $B_k=\bigcup_{i=1}^{k} \overline{t_i}$
  is a compact $3$--ball for all $k=1,\dots, 24n-1$. 

  Let $k=1,\dots, 24n-2$. If $\d t_{k+1}\cap B_k$ contains exactly $j$
  open $2$--simplices of $\T'$, then there is a unique open $(3-j)$-simplex
  $\sigma_k \subset \d t_{k+1}\setminus B_k$ of $\T'$. 
  There is a $(\T')^1$--Morse embedding $H\co S^2\times I\to S^3$
  with the following properties.
  \begin{enumerate}
  \item $H_0\subset t_1$.
  \item $H_{\frac{1}{24n}} = \d U(B_1)$ and $H$ has exactly four
    critical parameters in $[0,\frac{1}{24n}]$ 
    with respect to $(\T')^1$, with critical points
    in the vertices of $t_1$, and 
  \item $H_{\frac{k}{24n}}= \d U(B_k)$ for any $k=2,\dots 24n-1$, 
    and $H$ has exactly one critical parameter
    in $[\frac{k-1}{24n},\frac{k}{24n}]$ with respect to $(\T')^2$,
    namely with critical point in the barycenter of $\sigma_{k-1}$. Compare
    the $2$--dimensional sketch in Figure~\ref{fig:sweepshelling}.
    \begin{figure}
      \begin{center}
        {\psfrag{Bk-1}{$B_{k-1}$}\psfrag{sk-1}{$\sigma_{k-1}$}
            \psfrag{Hk-1}{$H_{\frac{k-1}{24n}}$}
          \psfrag{HK}{$H_{\frac{k}{24n}}$}\psfrag{t}{$t_k$}\psfrag{BK}{$B_k$}
          \epsfig{file=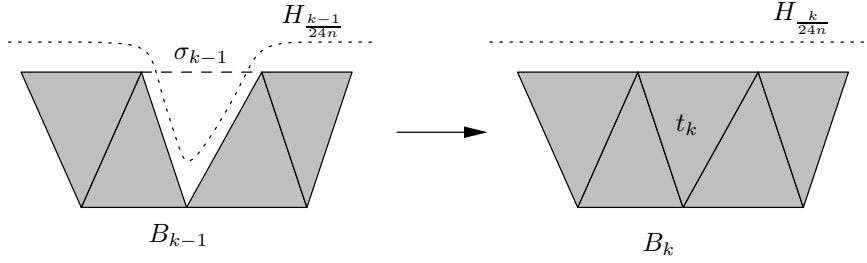}}
        \caption{Sweep-out along a shelling}
        \label{fig:sweepshelling}
      \end{center}
    \end{figure}
  \end{enumerate}
  By construction, $H$ has at most one critical
  point in each open simplex of $\T'$, namely in its barycenter.
  Thus the critical points of $H$ with respect to the $1$--skeleton of
  the dual cellular decomposition of $\T$ are its $n$ vertices
  and at most three critical points in each of its $2n$ open
  edges. This yields the lemma.    
\end{proof}

%****************************************************************

\subsection{The polytopality grows exponentially}  
\label{sec:expo}

This section is devoted to the proof of Theorem~\ref{thm:poly2}. 
We outline the construction after recalling a couple of notions.
A $2$--polyhedron $Q$ is \begriff{simple} or a \begriff{fake surface},
if the link of any point in $Q$ is homeomorphic to (i) a circle or
(ii) a circle with diameter or (iii) a complete graph with four
vertices.
A \begriff{$2$--stratum} of $Q$ is a connected component
of the union of points of type (i).
The points of type (iii) are the \begriff{intrinsic vertices} of $Q$.
A cellular decomposition $\C$ of a compact $3$--manifold is simple, 
if $|\C^2|$ is simple and $\C^0$ is a union of intrinsic vertices 
of $|\C^2|$. 
This corresponds to the classical notion of simple convex 
polytopes~\cite{ziegler}.

\begin{figure}
  \begin{center}
    \leavevmode
    {\psfrag{1}{$1$}\psfrag{2}{$2$}\psfrag{i+1}{$i+1$}\psfrag{i}{$i$}\psfrag{n}{$n$}\epsfig{file=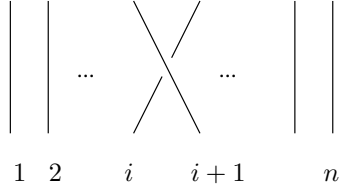}}
    \caption{The braid generator $\sigma_i$}
    \label{fig:ArtinGen}
  \end{center}
\end{figure}

For $n\in\mathbb N$, let $\mathcal B_n$ denote the group of braids with
$n$ strands; 
see~\cite{burde}, for instance. It
is generated by $\sigma_1,\dots,\sigma_{n-1}$, where $\sigma_i$
corresponds to a crossing of the $i$--th over\footnote{Note that in the
  literature appear different conventions on whether in $\sigma_i$ the
  $i$--th strand crosses \emph{over} or \emph{under} the $(i+1)$--th
  strand.}  
the $(i+1)$--th strand of the braid; see Figure~\ref{fig:ArtinGen}.

\begin{figure}
  \begin{center}
    \leavevmode
    {\psfrag{b}{$b$}\psfrag{b-}{$b^{-1}$}\psfrag{K1}{$K_1$}\psfrag{K2}{$K_2$}\epsfig{file=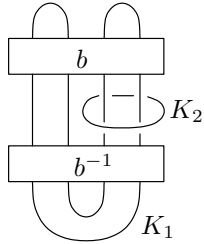}
    \caption{The link $K_1\cup K_2$}}
    \label{fig:closebraid}
  \end{center}
\end{figure}

Let $b=\sigma_{i_k}^{\epsilon_k} \cdot \sigma_{i_{k-1}}^{\epsilon_{k-1}} \cdots \sigma_{i_1}^{\epsilon_1}\in
\mathcal B_4$ be a braid with $k$ crossings, where $\epsilon_j\in
\{+1,-1\}$ for $j=1,\dots,k$. 
Let $K_1\cup K_2\subset S^3$ be the link that is defined in
Figure~\ref{fig:closebraid}. Both $K_1$ and $K_2$ are
unknots. 
We will use a result of Hass, Snoeyink and Thurston~\cite{hass},
providing a sequence of examples of a braid $b$ such that any spanning
disc for $K_1$ intersects $K_2$ in at least an exponential number of
points, in terms of $k$.
Let $V=S^3\setminus U(K_1)$, which is a solid torus containing $K_2$ as
a not necessarily trivial knot.  

\begin{figure}
    \begin{center}
      \leavevmode
      {\psfrag{K}{${K_1}$}\psfrag{z=1}{$z=1$}\psfrag{z=0}{$z=0$}\psfrag{x=i}{$x=i_j$}\psfrag{x=i+1}{$x=i_j+1$}\psfrag{y=j}{$y=j$}\psfrag{y=j+1}{$y=j+1$}\epsfig{file=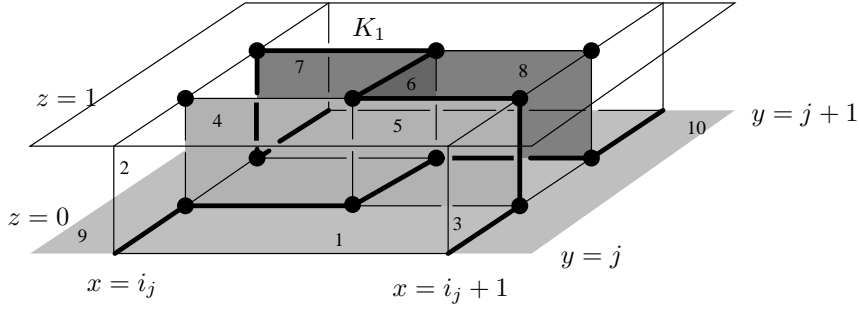}}
      \caption{Realization of the crossing $\sigma_{i_j}$ of $b$}
      \label{fig:crossings}
    \end{center}
  \end{figure}

\begin{figure}
    \begin{center}
      \leavevmode
      {\psfrag{K}{$K_1$}\psfrag{x=0}{$x=0$}\psfrag{x=5}{$x=5$}\psfrag{y=-k-2}{$y=-k-2$}\psfrag{y=-k-1}{$y=-k-1$}\epsfig{file=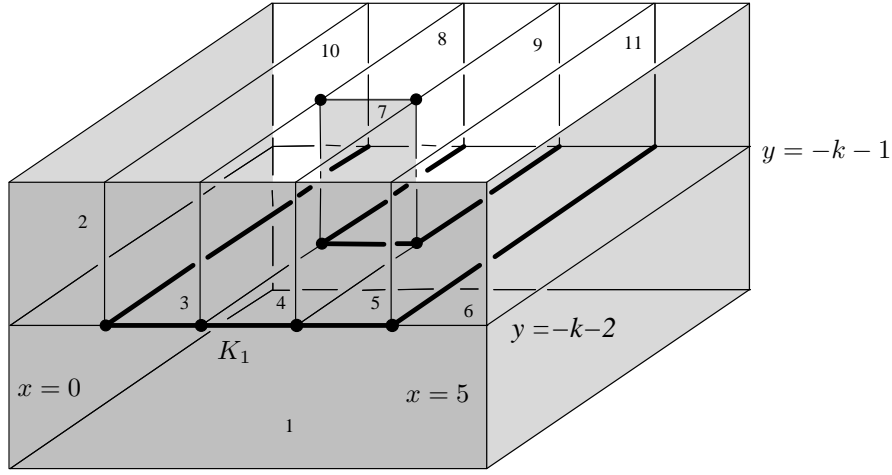}}
      \caption{Realization of the cups of $b$}
      \label{fig:cups}
    \end{center}
  \end{figure}

\begin{figure}
    \begin{center}
      {\psfrag{x=0}{$x=0$}\psfrag{K1}{$K_1$}\psfrag{K2}{$K_2$}\psfrag{x=5}{$x=5$}\psfrag{y=-1/2}{$y=-\frac 12$}\psfrag{y=1/2}{$y=\frac 12$}\epsfig{file=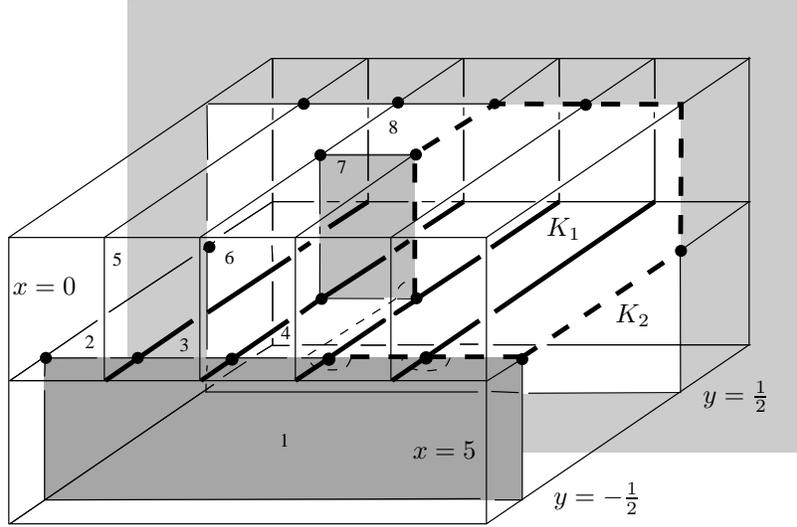}}
      \caption{Realization of $K_2$}
      \label{fig:middle}
    \end{center}
  \end{figure}

We explain the proof idea for Theorem~\ref{thm:poly2}.
We start with constructing a cellular
decomposition $\z_b$ of $V$ with $K_2\subset \z_b^1$. 
For the construction of $\z_b$, we put together the bricks shown in
Figures~\ref{fig:crossings}, \ref{fig:caps}, \ref{fig:cups}
and~\ref{fig:middle} (see Construction~\ref{cons:Qb}), and then drill
out a regular neighborhood of $K_1$ (see Lemma~\ref{lem:Q_bgut}).
Next we glue two modified copies of $\z_b$ together (see
Construction~\ref{con:Cb}) in order to
obtain a simple cellular decomposition of 
$S^3$ that is dual to a triangulation (see Lemma~\ref{lem:Z_bgut}). 
If one chooses $b$ according to~\cite{hass} then the polytopality of
the triangulation is ``very big'', yielding Theorem~\ref{thm:poly2}.
\begin{figure}
  \begin{center}
    \leavevmode
    \epsfig{file=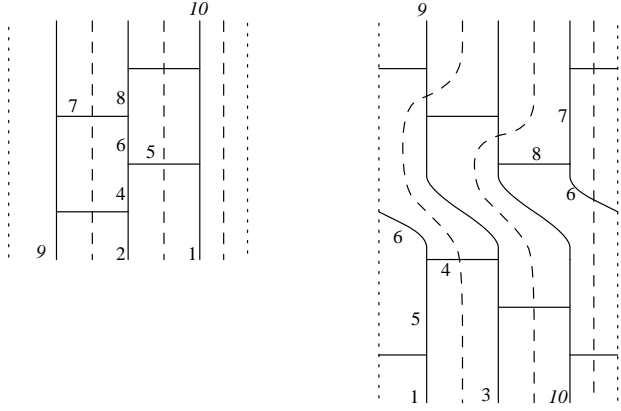}
    \caption{$\d V$ near a crossing of $b$}
    \label{fig:dcrossings}
  \end{center}
\end{figure}

\begin{construction}\label{cons:Qb}
  Let $(x,y,z)$ be coordinates for $\mathbb R^3=
  S^3\setminus\{\infty\}$. We define $W=\{0\le x\le 5, 
  -k-2\le y\le k+2, -1\le z\le 1\}$ and 
  \begin{eqnarray*}
   X &=& \{z=0\} \cup \{x\in \mathbb Z, z\ge 0\}\cup 
                      \{y=- \frac 12, z\le 0\} \cup \\ 
      && \{2\le x\le 3, y=0, z\ge 0\} \cup
             \{2\le x\le 3, y=-k - \frac 32, z\ge 0\}.
  \end{eqnarray*}        
  The simple $2$--polyhedron $P= \d W\cup (X\cap W) \cup (\{y=\frac 12\} 
  \setminus  W)$ has 28 intrinsic vertices. The unbounded $2$--stratum
  $\{y=\frac 12\}\setminus W$ is considered as a disc in $S^3$.    
  We define 
  \begin{eqnarray*} 
    P_b = &P&  \cup  \bigcup_{j=1}^k
           \Big(\{i_j\le x\le i_j+1, \pm y=j + \frac 13, 0\le z\le 1\}
           \cup\\  
         &&\quad \{i_j\le x\le i_j+1, \pm y=j + \frac 23, 0\le z\le 1\}
         \cup\\ 
         && \quad\{x=i_j+\frac 12, j+\frac 13 \le \pm y\le j+\frac 23, 
               0\le z\le 1\} \Big).
  \end{eqnarray*}
\begin{figure}
    \begin{center}
      \leavevmode
      {\psfrag{K}{$K_1$}\psfrag{x=0}{$x=0$}\psfrag{x=5}{$x=5$}\psfrag{y=k+1}{$y=k+1$}\psfrag{y=k+2}{$y=k+2$}\epsfig{file=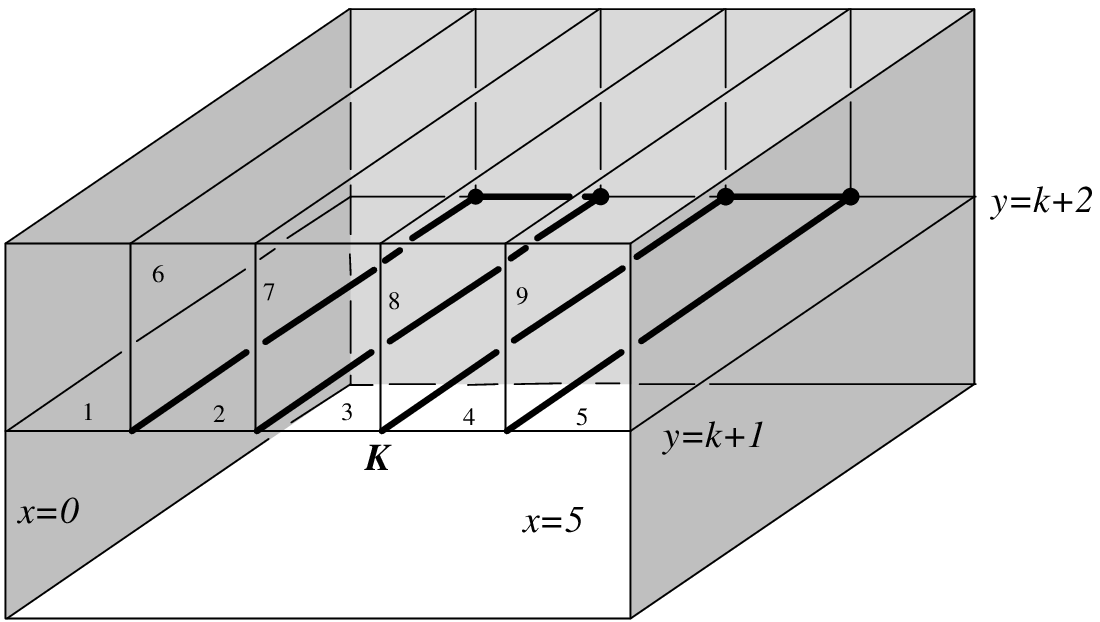}}
      \caption{Realization of the caps of $b$}
      \label{fig:caps}
    \end{center}
  \end{figure}
\begin{figure}
  \begin{center}
    \leavevmode
    \epsfig{file=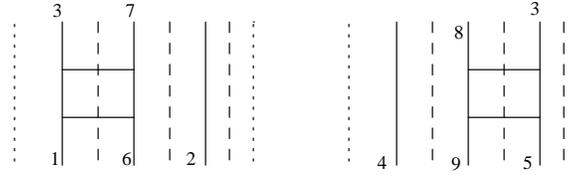}
    \caption{$\d V$ in the caps of $b$}
    \label{fig:dcaps}
  \end{center}
\end{figure}

  One observes that $P_b$ is the $2$--skeleton of a simple cellular
  decomposition of $S^3$ with $24k+28$ vertices that is dual to a
  triangulation. Any 
  crossing of the braid $bb^{-1}$ is realized in $P_b^1$; see
  Figure~\ref{fig:crossings}. The figure shows the crossing
  $\sigma^{+1}_{i_j}$ of $b$, where $b$ is bold. 
  Here and in all
  subsequent figures, thick dots indicate intrinsic vertices of
  simple $2$--polyhedra. %
  Also the ``caps'' and
  ``cups'' that yield $K_1$ from $bb^{-1}$ are realized in $P_b^1$; see
  Figures~\ref{fig:caps} and~\ref{fig:cups}. Thus we can assume
  $K_1\subset P_b^1$.\qed
\end{construction}
\begin{lemma}\label{lem:Q_bgut}
  The $2$--polyhedron $Q_b = (P_b\setminus U(K_1))\cup \d U(K_1)\subset
  V$ is the $2$--skeleton of a simple cellular decomposition $\z_b$ of $V$
  with $< 48 k + 56$ vertices that is dual to a triangulation of
  $V$. A representative of $K_2\subset V$ is formed by $11$ edges of
  $\z_b$. 
\end{lemma}
\begin{figure}
  \begin{center}
    \leavevmode
    \epsfig{file=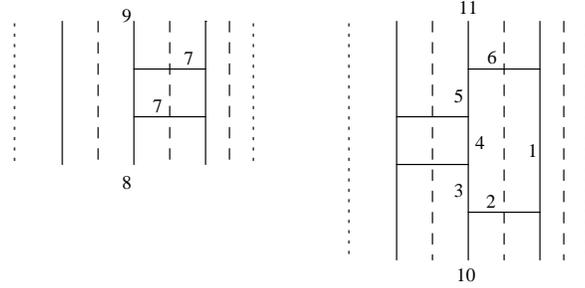}
    \caption{$\d V$ in the cups of $b$}
    \label{fig:dcups}
  \end{center}
\end{figure}
\begin{proof}  
  Since $P_b$ is the $2$--skeleton of a cellular decomposition $S^3$ dual
  to a triangulation, any connected component of $V\setminus Q_b$ is a
  ball, the intersection of any two of these compact balls is
  connected, and the $2$--strata of $Q_b$ in the
  interior of $V$ are discs.

  We show that the closure of any $2$--stratum of $Q_b$ in $\d V$ is a
  disc, and the intersection of any two of these discs is connected.
  A $2$--stratum of $Q_b$ in $\d V$ is a connected component of $\d U(K_1)
  \setminus P_b$.  Since $K_1$ is not contained in the boundary of a
  single connected component of $S^3\setminus P_b$, the closure of any
  component of $\d U(K_1) \setminus P_b$ is a disc.
  Since $K_1$ is not a union of two arcs that are each contained in the
  boundary of a connected component of $S^3\setminus P_b$, the
  intersection of the closures of two connected components of $\d U(K_1)
  \setminus P_b$ is connected.

  Therefore $Q_b$ is the $2$--skeleton (including $\d
  V$) of a simple cellular decomposition $\mathcal Z_b$ of $V$, and the
  dual of  $\mathcal Z_b$ has no loops or multiple edges, thus,  is
  a triangulation.  
  Any vertex of $P_b$ in $K_1$ gives rise to two vertices of $\mathcal
  Z_b$. Thus  $\mathcal Z_b$ has $< 48 k + 56$ vertices. 
  The $1$--skeleton of $\z_b$ contains $K_2$ as a path of $11$ edges, namely
  seven edges corresponding to edges of $P_b$ (the thick dotted lines in
  Figure~\ref{fig:middle}) and four edges in $P_b\cap \d U(K_1)\subset
  \z_b^1$ (the thin dotted lines in the figure). 
\end{proof}

Our plan is to glue two copies of $\mathcal Z_b$ together, to obtain a
simple cellular decomposition of $S^3$. In order to keep a bound for
the number of vertices that is linear in $k$, we need an ``adaptor''
between the two copies of $\mathcal Z_b$.
The construction of the adaptor is based on the following lemma.

\begin{lemma}\label{lem:arrangeGamma}
  There is a disjoint union $\Lambda\subset \d V$ of three meridians
  of $V$ such that $\#(\Lambda \cap \mathcal Z_b^1) = 22 k + 16$ and
  the intersection of any component of $\d V\setminus \Lambda$ with
  any $2$--cell of $\z_b$ in $\d V$ is connected.
\end{lemma}
%**************
\begin{figure}
  \begin{center}
    \epsfig{file=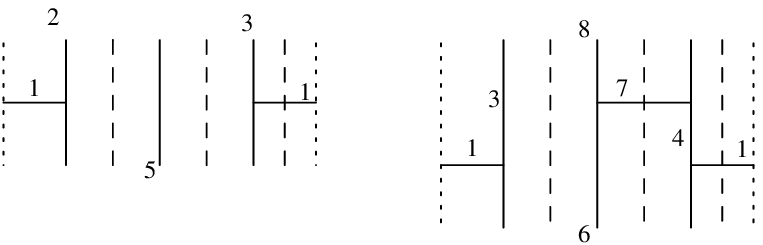}
    \caption{$\d V$ in $\{-1 < y < 1\}$}
    \label{fig:dmiddle}
  \end{center}
\end{figure}

\begin{proof}
  We construct $\Lambda\subset \d V$ according to
  Figures~\ref{fig:dcrossings}, \ref{fig:dcaps}, \ref{fig:dcups}
  and~\ref{fig:dmiddle}.
  Figure~\ref{fig:dcrossings} shows the two annuli contained in $\d V$
  that correspond to the two sub-arcs of $K_1$ shown in
  Figure~\ref{fig:crossings}; the annuli are cut along the dotted lines
  (left and right side of the rectangles in the 
  figure). 
  The figure shows the pattern of $\d V\cap \mathcal Z_b^1$, where the
  numbers in Figure~\ref{fig:dcrossings} at the edges of 
  $\mathcal  Z_b$ in $\d V$ correspond to the numbers in
  Figure~\ref{fig:crossings} at the $2$--strata of $P_b$. 
  The broken lines indicate $\Lambda$. One 
  sees $11$ points of $\Lambda\cap\mathcal Z_b^1$. Thus the $2k$
  crossings of $bb^{-1}$ contribute to $22 k$ points in $\Lambda\cap
  \mathcal Z_b^1$. 

  Similarly, Figures~\ref{fig:dcaps} and~\ref{fig:dcups} show the
  parts of $\d V$ corresponding to the sub-arcs of $K_1$ in
  Figures~\ref{fig:caps} and~\ref{fig:cups}. We see $4$ respectively $6$
  points of $\Lambda \cap \mathcal Z_b^1$. In
  Figure~\ref{fig:dmiddle}, we show two of the four parts of $\d V$
  corresponding to Figure~\ref{fig:middle}, and obtain by symmetry 6
  points of $\Lambda\cap \mathcal Z_b^1$.

  The broken lines in the figures close up to three meridians (forming $\Lambda$), since
  the writhing number of $K_1$ vanishes. We have $\#(\Lambda\cap
  \z_b^1) = 22k + 4 + 6 + 6 = 22k+ 16$.
  The second claim of the lemma follows, since any $2$--cell of $\z_b$ in
  $\d V$ meets $\Lambda$ in at most two arcs.  
\end{proof}

\begin{figure}
  \begin{center}
    {\psfrag{z}{$\z_b$}\psfrag{z'}{$\z_b'$}\psfrag{V}{$\d V$}\epsfig{file=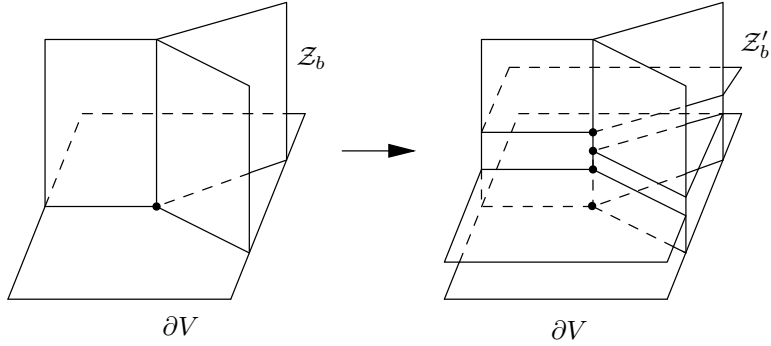}}
    \caption{Subdividing $3$--cells of $\z_b$}
    \label{fig:subdivide}
  \end{center}
\end{figure}

\begin{construction}\label{con:Cb}
  We subdivide each $3$--cell of $\z_b$ that meets $\d V$ by inserting a
  copy of all $2$--cells of $\z_b$ in $\d V$, where the copies are
  pairwise disjoint\footnote{This corresponds to stellar subdivisions
    along some edges of the dual triangulation of $\z_b$.}; see
  Figure~\ref{fig:subdivide}.  By Lemma~\ref{lem:Q_bgut}, we obtain a
  simple cellular decomposition $\z'_b$ of $V$ with $< 4\cdot(48k+56)
  = 192 k + 224$ vertices dual to a triangulation such that the
  intersection of any compact $3$--cell with $\d V$ is connected. A
  representative of $K_2$ is formed by a path of at most $4\cdot 11$
  edges of $\z'_b$.

  Decompose $S^3= V_0\cup_{\d V_0} (S^1\times S^1\times I) \cup_{\d V_1}
  V_1$, where $V_0,V_1$ are solid tori and $S^1\times \{*\}\times \{0\}$
  (resp. $\{*\}\times S^1\times \{1\}$) is a meridian for $V_0$
  (resp. $V_1$). 
  For $i=0,1$, let $\mathcal Z_{b,i}$ be a cellular
  decomposition of $V_i$ isomorphic to  $\mathcal Z'_b$. 

  Let $\Gamma\subset S^1\times S^1$ be a union of three copies of
  $S^1\times \{*\}$ and $\{*\}\times S^1$, meeting in nine points.
  We attach $\Gamma\times \{0\}$ and $\Gamma\times \{1\}$ at $\d V_0$
  and $\d V_1$ as follows.
  Choose the three meridians in $\Gamma\times \{0\}\subset \d V_0$ for
  $V_0$ and the three meridians in $\Gamma\times \{1\}\subset \d V_1$
  for $V_1$ according to Lemma~\ref{lem:arrangeGamma}. 
  We choose the three longitudes in $\Gamma\times \{0\}\subset \d V_0$
  for $V_0$ (resp.\ in $\Gamma\times \{1\}\subset \d V_1$ for $V_1$)
  so that each of them intersects $\z_{b,0}^1$ (resp. $\z_{b,1}^1$) in
  three points, and $\Gamma\times \{0\}$ (resp. $\Gamma\times \{1\}$)
  meets each $2$--cell of $\z_{b,0}^1$ (resp.\ of $\z_{b,1}^1$) in at
  most two arcs.
  
  We define $C = \mathcal Z^2_{b,0}\cup (\Gamma\times I) \cup
  \mathcal Z^2_{b,1}\subset S^3$. 
  By Lemma~\ref{lem:arrangeGamma} and by construction of $\Gamma$, $Z$
  has $< 2\cdot (192 k + 224 + 22 k + 16 + 18) = 428 k +  516$ 
  intrinsic vertices.
  It has nine $4$--valent edges, corresponding to the vertices of
  $\Gamma$. 
  We perturb $C$ along these edges, which increases the number of 
  vertices by 18, and obtain a simple $2$--polyhedron $C'$.\qed
\end{construction}

\begin{lemma}\label{lem:Z_bgut}
  The $2$--polyhedron $C'$ is the $2$--skeleton of a {simple} cellular
  decomposition $\C_b$ of $S^3$ that is dual to a triangulation with
  $< 428 k + 534$ tetrahedra.
\end{lemma}
\begin{proof} 
  By Lemmas~\ref{lem:Q_bgut} and~\ref{lem:arrangeGamma}, the
  closure of each $2$--stratum of $Z'$ is a disc and 
  the closure of any connected component of $S^3\setminus C'$ is a
  ball. Thus $Z'$ is the $2$--skeleton of a simple cellular decomposition
  of $S^3$.

  Let $X_1,X_2$ be the closures of two connected components of
  $S^3\setminus C'$. It remains to show that $X_1\cap X_2$ is connected.
  If $X_1\subset V_1$ and $X_2\subset V_2$ then $X_1\cap X_2 = \emptyset$.
  If $X_1\cup X_2\subset S^1\times S^1\times I$ then $X_1\cap X_2$ is
  connected by construction of $\Gamma$.
  According to Lemma~\ref{lem:Q_bgut}, if $X_1\cup X_2\subset V_i$ for
  $i=1,2$ then $X_1\cap X_2$ is connected.
  Finally, let $X_1\subset V_i$ and $X_2\subset S^1\times S^1\times I$. 
  Then $X_1\cap X_2$ is connected by
  Lemma~\ref{lem:arrangeGamma} and since we have subdivided  
  the $3$--cells of $\z_b$ that meet $\d V$.    
\end{proof}

\medskip
\noindent
\emph{Proof of  Theorem~\ref{thm:poly2}.} 
Let $m\in \mathbb N$.
We set $b=(\sigma_1\sigma_2^{-1})^m \in \mathcal B_4$ in order to
apply the results of~\cite{hass}. 
Let $\C_m = \C_{b}$ be as in Lemma~\ref{lem:Z_bgut};
it is dual to a triangulation $\T_m$ of $S^3$ with $\le  428 k + 534 =
856 m +534$ tetrahedra.

It remains to show that $p(\T_m)> 2^{m-1}$.
For $i=0,1$, let $K_{2,i}\subset \C_m^1$ be the copy of $K_2$ in
$V_i$. 
Let $L_m = K_{2,0}\cup K_{2,1}$, which by Construction~\ref{con:Cb} is
a link formed by at most 88 edges of $\C_m$.
Let $H\co S^2\times I\to S^3$ be a $\C_m^1$--Morse embedding with
$\C_m^2\subset H(S^2\times I)$.
There is a parameter $\xi\in I$ for which $H_\xi$ contains a meridional disc for $V_0$ or for $V_1$.
By~\cite{hass}, any meridional disc in $V_i$ intersects 
$K_{2,i}$  in $\ge 2^{m-1}$ points. 
Since $L_m \subset \C_m^1$ and since $\#(H_\xi\cap
L_m)$ changes by $\pm 2$ at any critical parameter, we have $p(\T_m)>2^{m-1}$, which proves
Theorem~\ref{thm:poly2}.

%***********************************************************************
%***********************************************************************

\section{Contractions and expansions}
\label{sec:hauptbeweise}

This section is devoted to the proofs of Theorems~\ref{thm:basic3}
and~\ref{thm:basic2}. The proofs are outlined in the introduction.

\subsection{Canceling pairs of critical parameters} 
\label{sec:minmax}

The setting of this subsection is more general as actually needed to
prove Theorems~\ref{thm:basic3} and~\ref{thm:basic2}.
Let $M$ be a closed orientable $3$--manifold with a cellular
decomposition $\C$ that is dual to a triangulation of $M$. Let $S$ be
a closed surface.  The aim of this section is to deform a
$\C^1$--Morse embedding $H\co S\times I \to M$ by isotopy so that
$c(H,\C^1)$ is unchanged and $c(H,\C^2)$ is bounded in terms of
$c(H,\C^1)$.

\begin{figure}
  \begin{center}
    {\psfrag{E}{$E^+$}\psfrag{V}{$V^+$}\epsfig{file=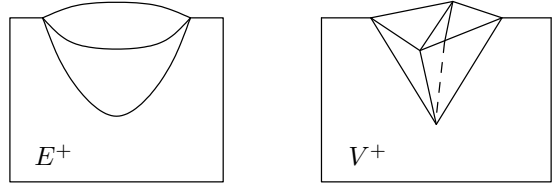}}
    \caption{Types of critical points in $\C^1$}
    \label{fig:crit1}
  \end{center}
\end{figure}

\begin{figure}
  \begin{center}
    {\psfrag{F+}{$F^+$}\psfrag{F0}{$F^0$}\epsfig{file=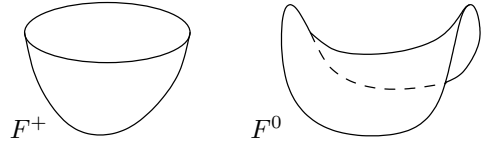}}
    \caption{Types of critical points in a $2$--cell}
    \label{fig:crit2}
  \end{center}
\end{figure}

A $\C^1$--Morse is essentially the same as an embedding of the product
$S\times I$ in general position with respect to $\C^1$. We extend this
to the notion of $\C^2$--Morse embeddings. The critical points in an
open $2$--cell are of the usual Morse types (local minimum, local
maximum, or saddle point). Since general position with respect to
$\C^2$ is not sufficient for our purposes, we further restrict the
occuring types of critical points in $\C^1$ as in the following two
definitions. Compare Figures~\ref{fig:crit1} and~\ref{fig:crit2}.

\begin{definition}\label{def:types}
  Let $H\co S\times I \to M$ be a $\C^1$--Morse embedding.
  Let $p_0$ be a critical point of $H$ with respect to $\C^2$,
  associated to a critical parameter $\xi_0\in I$. If there are local
  coordinates $(x,y,z)$ around $p_0$
  such that $H_{\xi}\cap U(p_0) = \{z=\xi-\xi_0\}\cap U(p_0)$ for
  $\xi\in I$, and
  $\C^2\cap U(p_0)$ equals
  \begin{itemize}
  \item $\{z= x^2 + y^2\}\cap U(p_0)$, then $p_0$ (resp.\
    $\xi_0$) is of type $F^+$.
  \item $\{-z= x^2 + y^2\}\cap U(p_0)$, then $p_0$ (resp.\
    $\xi_0$) is of type $F^-$.
  \item $\{z= x^2 - y^2\}\cap U(p_0)$, then $p_0$ (resp.\
    $\xi_0$) is of type
    $F^0$.
  \item $\left(\{z= x^2 + |y|\}\cup \{z\le  x^2,
    y=0\}\right)\cap U(p_0)$, then $p_0$ (resp.\
    $\xi_0$) is of type $E^{+}$.
  \item $\left(\{-z= x^2 + |y|\}\cup \{-z\le  x^2,
    y=0\}\right)\cap U(p_0)$, then $p_0$ (resp.\
    $\xi_0$) is of type $E^{-}$.
  \item $\left(\{z=|x|+|y|\} \cup \{z\le |x|, y=0\} \cup \{z\ge |y|,
      x=0\}\right)\cap U(p_0)$, then $p_0$ (resp.\ $\xi_0$) is of type
    $V^+$. 
  \item $\left(\{-z=|x|+|y|\} \cup \{-z\le |x|, y=0\} \cup \{-z\ge |y|,
      x=0\}\right)\cap U(p_0)$, then $p_0$ (resp.\ $\xi_0$) is of type
    $V^-$. 
  \end{itemize}
\end{definition}

\begin{definition}
  Let $H\co S\times I\to M$ be a $\C^1$--Morse embedding. If any critical
  point of $H$ with respect to $\C^2$ is of type $F^\pm$,
  $F^0$,$E^\pm$ or $V^\pm$, then $H$ is a \begriff{$\C^2$--Morse
    embedding}. 
\end{definition}
Note that for each critical point $p_0\in \C^2$ of $H$ corresponding
to a critical parameter $\xi_0$, there is a \emph{unique} open $2$--cell
$c$ of $\C$ such that $U(p_0)\cap H_{\xi_0} \cap c \not= \emptyset$.
We will use this fact in Construction~\ref{con:alpha} and
Lemma~\ref{lem:minweg}.

Let $H\co S\times I\to M$ be a $\C^2$--Morse embedding such that $H_\xi$
splits $M$ into two pieces for all $\xi\in I$.
Let $B^+,B^-$ be the two components of $M\setminus
H(S\times I)$, with $\d B^+ = H_1$ and $\d B^- = H_0$. For $\xi\in I$,
let $B^+(\xi)$ (resp.\ $B^-(\xi)$) be the closure of the component of
$M\setminus H_\xi$ that contains $B^+$ (resp.\ $B^-$).

Our \textbf{general hypothesis} in this section is that both $0$ and $1$
are non-critical parameters of $H$ with respect to $\C^2$, that
$H_0\cap \C^2\not=\emptyset\not=H_1\cap \C^2$, and that, for any
open $3$--cell $X$ of $\C$, any circle in $H_0\cap \d X$ (resp. in
$H_1\cap \d X$) bounding a disc in $B^+(0)\cap \d X$ (resp. in
$B^-(1)\cap \d X$) also bounds a disc in $H_0\cap X$ (resp. in
$H_1\cap X$).

We will change $H$ by isotopy into a
$\C^2$--Morse embedding $\tilde H\co S\times I\to M$ so that
$c(\tilde H,\C^2)$ is subject to an upper bound in terms of
$c(H,\C^1)$.  
Whenever $H$ has a critical point of type $F^\pm$, then it ``cancels''
with a critical point of type $F^0$ (see Lemma~\ref{lem:splitdisc}). 
We change $H$ according to Figures~\ref{fig:moveH}
and~\ref{fig:moveH2}, removing the canceling pair of critical points (see
Lemma~\ref{lem:minweg}).  

Let $p_0\in \C^0$ be a critical point of type $F^+$ that belongs to a
critical parameter $\xi_0$ of $H$.  There is a unique open $3$--cell $X$
of $\C$ such that $U(p_0)\cap H_{\xi_0} \cap X \not=\emptyset$.
\begin{lemma} \label{lem:splitdisc}
\begin{enumerate}
\item  There is a nonempty open interval $(\xi_0, \overline{\xi_0})$
  such that for any $\xi \in (\xi_0, \overline{\xi_0})$ there is
  a connected component $\gamma(\xi)$ of $H_\xi\cap \d X$ that is a
  circle, such that $\gamma(\xi)$  varies continously in $\xi$ and
  $\lim_{\xi \to \xi_0} \gamma(\xi) = p_0$. 

\item  If $\overline{\xi_0}$ is maximal then $\gamma(\overline{\xi_0}) =
  \lim_{\xi \nearrow \overline{\xi_0}} \gamma(\xi)$ is a circle that is not a
  connected component of $H_{\overline{\xi_0}}\cap \d X$ and contains a
  critical point $\overline{p_0}$ of $H$ of type $F^0$.   
\end{enumerate}
\end{lemma}

\begin{proof}
  For the first part of the lemma, let $\epsilon > 0$ be small enough
  such that there is no critical parameter of $H$ with respect to
  $\C^2$ in $(\xi_0,\xi_0+\epsilon)$. Then, $H_{\xi_0+\epsilon} \cap
  \d X$ contains a small circle $\gamma(\xi_0+\epsilon)$ around $p_0$.
  For $\xi\in (\xi_0,\xi_0+\epsilon)$, one obtains $\gamma(\xi)$ from
  $\gamma(\xi_0+\epsilon)$ by isotopy induced by $H$.
  
  For the second part of the lemma, note that $\gamma(\xi_0+\epsilon)$
  bounds a disc in $\d X\cap B^-(H_{\xi_0+\epsilon})$. By induction,
  $\gamma(\overline{\xi_0})$ bounds a disc in $\d X\cap
  B^-(H_{\overline{\xi_0}})$ and does not bound a disc in
  $H_{\overline{\xi_0}}\cap X$. Thus by the general hypothesis,
  $\overline{\xi_0}<1$.
  By maximality of $\overline{\xi_0}$, either $\gamma(\overline{\xi_0})$
  is a critical point of $H$ of type $V^-$, $E^-$ or $F^-$, or
  $\gamma(\overline{\xi_0})$ is a circle 
  containing a critical point of $H$ of type $F^0$. 
  Since $H_0\cap \C^2\not= \emptyset$ by our general hypothesis, it
  follows by induction that the
  connected component of $H_{\xi}\cap X$ containing $\gamma(\xi)$ in its
  boundary is not a disc, for $\xi \in (\xi_0, \overline{\xi_0})$. 
  Thus the critical point  in
  $\gamma(\overline{\xi_0})$ is of type $F^0$.  
\end{proof}
We keep the notations of the preceding lemma and assume that
$\overline{\xi_0}$ is maximal. Let $\xi_0<\xi_1 <\dots
<\xi_k=\overline{\xi_0}$ be the critical parameters of $H$ with
respect to $\C^2$ in the closed interval $[\xi_0,\overline{\xi_0}]$,
corresponding to the critical points $p_0,p_1,\dots,p_k$. The
following Construction~\ref{con:alpha} yields two compact arcs
$\alpha_1,\alpha_2\subset \d X$ so that for $\xi \in [\xi_0,\xi_k]$
and $m=1,2$ holds
\begin{enumerate}
\item $\alpha_m\cap H_\xi$ is empty or a single point
$\alpha_m(\xi)$,
\item $\alpha_1\cup \alpha_2$ contains at
  most one critical point of $H$ of type $F^+$,
\item $\alpha_m\cap \C^1$ consists of critical points of $H$ of type $E^+$ and
  $V^+$,  and 
\item $\alpha_2(\xi)\in \gamma(\xi)$, and either $\alpha_1(\xi)\not\in \gamma(\xi)$ or $\alpha_1(\xi)=\alpha_2(\xi)=p_k$.  
\end{enumerate}

\begin{construction} \label{con:alpha}
  We construct $\alpha_1,\alpha_2$ iteratively.  Let $\alpha_1\cap
  B^+(\overline{\xi_0})=\alpha_2 \cap B^+(\overline{\xi_0}) = p_k$.
  Let $i\in \{0,1,\dots,k\}$. For both $m=1,2$, suppose that
  $\alpha_m\cap H(S^2\times [\xi_i,\xi_k])$ is already constructed.
  If $\alpha_m(\xi_i)$ is a critical point of $H$ of type $F^+$, then
  we stop the construction and define $\xi_0'=\xi_i$. 
  Otherwise, by construction either $\alpha_m(\xi_i)\not\in \C^1$ or
  $\alpha_m(\xi_i)$ is a critical point of $H$ of type $E^+$ or $V^+$.
  Thus by definition of $\C^2$--Morse embeddings, there is
  a unique open $2$--cell $c_m\subset \d X$ of $\C$ with
  $U(\alpha_m(\xi_i))\cap H_{\xi_i}\cap c_m \not=\emptyset$.  
  We extend $\alpha_m$ by an arc in $\overline{c_m} \cap
  H\left(S\times [\xi_{i-1},\xi_i]\right)$, so that 
  \begin{enumerate}
  \item $\alpha_1(\xi)\not\in \gamma(\xi)$ and $\alpha_2(\xi)\in\gamma(\xi)$ for
    $\xi\in (\xi_{i-1},\xi_i)$,
  \item if $\alpha_m(\xi_{i-1}) \in \C^1$ then
    $\alpha_m(\xi_{i-1})=p_{i-1}$, and 
  \item $\alpha_m(\xi_{i-1}) = p_{i-1}$ if and only if
    $U(\alpha_m(\xi_{i-1})) \cap  H_{\xi_{i-1}} \cap c_m=\emptyset$.
  \end{enumerate}
  This is possible, since if $\alpha_2(\xi_i)\in \gamma(\xi_i)$ then
  we can stay in $\gamma(\xi)$, if $\alpha_m(\xi_{i-1})$ is not a
  critical point then we can avoid to run into $\C^1$, and if
  $U(p_{i-1}) \cap H_{\xi_{i-1}} \cap c_m \not=\emptyset$ then we can
  avoid to run into the critical point $p_{i-1}$.
  It follows by property~(c) that $\alpha(\xi_{i-1})=p_{i-1}$
  only if $p_{i-1}$ is of type $F^+$, $E^+$ or $V^+$.\qed
\end{construction}

\begin{lemma}\label{lem:minweg}
  There is a $\C^2$-Morse embedding $\widetilde H\co  S\times I\nach M$
  isotopic to $H$ without critical points of type $F^\pm$ such that
  $\widetilde H(S\times \{0,1\}) = H(S\times \{0,1\})$ and
  $c(\widetilde H,\C^1) = c(H,\C^1)$.
\end{lemma}
\begin{proof}
  Let $\xi_0\in I$ be a critical parameter of $H$ of type $F^+$, and let
  $p_0$ be the corresponding critical point. 
  We use the notations of Construction~\ref{con:alpha}.
  If $\xi_0<\xi_0'$ then we replace $p_0$ by the critical point of $H$
  of type $F^+$ that corresponds to the critical parameter
  $\xi_0'$. Since there are only finitely many critical parameters in
  $[\xi_0,\xi_0']$, we can choose $p_0$ so that
  $\xi_0=\xi_0'$, by repeating this process. Thus we can assume that
  $\alpha_2(\xi_0) = p_0$. 

  Let $(x,y,z)$ be local coordinates around $p_k$ as in
  Definition~\ref{def:types}. 
  Let $\epsilon>0$ and $U(p_k)$ such that
  $H_{\xi_k-\epsilon}\cap U(p_k)$ is a disc, and
  $\xi_0$ (resp. $\xi_k$) is the only critical parameter of 
  $H$ with respect to $\C^2$ in $[\xi_0-\epsilon,\xi_0+\epsilon]$
  (resp.\ $[\xi_k-\epsilon,\xi_k+\epsilon]$). 
  Define $D' = \{-\epsilon \le z \le - y^2, x=0\}\subset U(p_k)$. We
  have $\d D'\subset H_{\xi_k-\epsilon}\cup \C^2$.  
  There is a compact disc $D\subset X\cap B^-(H_{\xi_k-\epsilon})$ such that 
  $D\cap H_\xi$ is a copy
  of $\gamma(\xi)$, for all $\xi\in (\xi_0,\xi_k-\epsilon]$.

  Let $\alpha_1,\alpha_2\subset \d X$ be as in
  Construction~\ref{con:alpha}.
  We change $H$ into an embedding $\widetilde H\co S\times I \to M$ in
  the following way, compare Figures~\ref{fig:moveH}
  and~\ref{fig:moveH2}.   
  \begin{figure}
    \begin{center}
      \leavevmode
      {\psfrag{H0}{$H_{\xi_0}$}\psfrag{H+}{$H_{\xi_k+\epsilon}$}\psfrag{H-}{$H_{\xi_k-\epsilon}$}\psfrag{pk}{$p_k$}\psfrag{p0}{$p_0$}\psfrag{a1}{$\alpha_1$}\psfrag{a2}{$\alpha_2$}\psfrag{D}{$D$}\epsfig{file=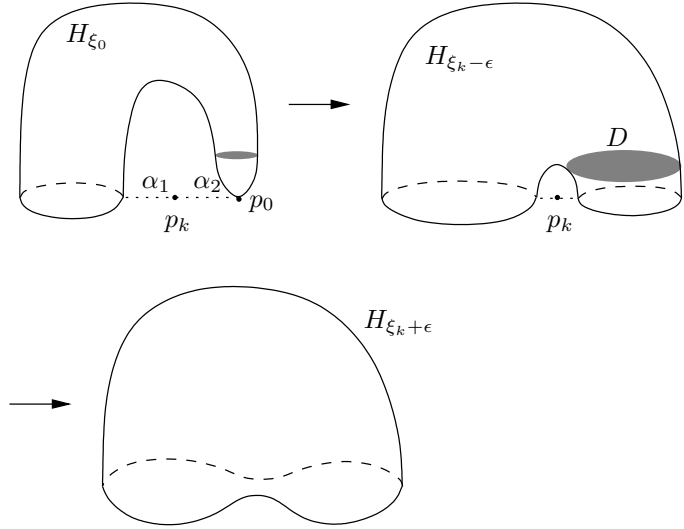}}
      \caption{The moves of $H$}
      \label{fig:moveH}
    \end{center}
  \end{figure}
  \begin{figure}
    \begin{center}
      \leavevmode
      {\psfrag{0}{$\widetilde H_{\xi_0}$}\psfrag{+}{$\widetilde H_{\xi_k+\epsilon}$}\psfrag{-}{$\widetilde H_{\xi_k-\epsilon}$}\epsfig{file=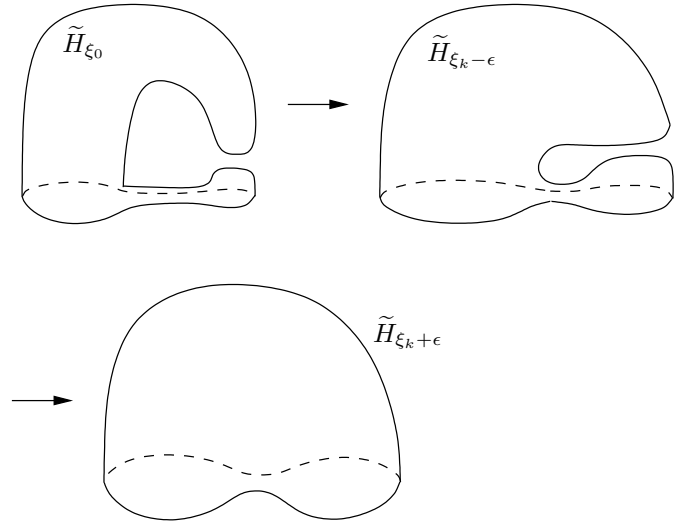}}
      \caption{The moves of $\widetilde H$}
      \label{fig:moveH2}
    \end{center}
  \end{figure}

  \begin{enumerate}
  \item For $\xi \in [0,\xi_0-\epsilon]$, let $\widetilde
    H(\cdot,\xi)\equiv H(\cdot,\xi)$. 
  \item When the parameter $\xi$ increases from $\xi_0-\epsilon$ to
    $\xi_0$, $\widetilde H$ pushes a finger 
    along $\alpha_1\cup\alpha_2$ towards $p_0$, so that 
    $\widetilde H_{\xi_0}
      = \d\left( 
          \left(B^-({\xi_0})\setminus U(D)\right)
          \cup U(\alpha_1\cup \alpha_2)
      \right).$
  \item For $\xi\in [\xi_0, \xi_k-\epsilon]$, set
     $\widetilde H_\xi = 
        \d\left( 
          \left(B^-(\xi)\setminus U(D)\right) 
          \cup U(\alpha_1\cup \alpha_2) 
        \right).$
  \item In $[\xi_k-\epsilon,\xi_k+\epsilon]$,
    $\widetilde H$ induces an isotopy mod $\C^2$ with support
    in $U(D'\cup D)$, relating  $\widetilde  
    H_{\xi_k-\epsilon}$ with $\widetilde
    H_{\xi_k+\epsilon} = H_{\xi_k+\epsilon}$. 
  \item For $\xi \in [\xi_k+\epsilon,1]$, let $\widetilde
    H(\cdot,\xi)\equiv H(\cdot,\xi)$. 
  \end{enumerate}

  When we push a finger along $\alpha_1\cup\alpha_2$, critical points of
  $\widetilde H$ of type $E^+$ or $V^+$ occur at
  $(\alpha_1\cup\alpha_2)\cap \C^1$.  
  By Construction~\ref{con:alpha}, these are also critical points of $H$
  of type $E^+$ and $V^+$. 
  Critical points of $\widetilde H$ in $\C^1\cap \d
  U(\alpha_1\cup\alpha_2)$ do not occur\footnote{This would be the case
    if $(\alpha_1\cup\alpha_2)\cap \C^1$ would contain non-critical
    points or critical points of type $E^-$ or $V^-$.}.
  It follows that the critical points of $H$
  and $\widetilde H$ in $\C^1$ coincide, $c(\widetilde H,\C^1) = c(H,\C^1)$,
  although the order of the corresponding critical \emph{parameters}
  has changed.  

  One sees that $\widetilde H$ has exactly two critical points with
  respect to $\C^2$ less than $H$, namely $p_0$ and $p_k$. One iterates
  this construction and removes all critical points
  of type $F^+$ of $H$. By the symmetric construction, one also removes
  all critical points of type $F^-$ of $H$, and the lemma follows. 
\end{proof}
 
After getting rid of the critical points of type $F^\pm$, it remains
to estimate the number of critical points of type $F^0$.
\begin{lemma}\label{lem:saddlebound}
  Assume that $H$ has no critical points of type $F^\pm$. Then $H$ has
   $\le \chi(B^-\cap \C^2) + \chi(B^+\cap \C^2) 
       - \chi(\C^2) + \chi(\C^0) + c(H,\C^1)$ critical
   points of type $F^0$.  
\end{lemma}
\begin{proof}
  For any $\xi\in I$, define
  $P_\xi = \left(\C^2\cap B^+(\xi)\right)\cup H_\xi$. The
  homeomorphism type of $P_\xi$ changes only at critical parameters of
  $H$ with respect to $\C^2$. 
  Let $\xi_0$ be a critical parameter of $H$ with respect to $\C^2$.
  We choose $\epsilon>0$ so that $\xi_0$ is the only critical
  parameter of $H$ in the interval 
  $[\xi_0-\epsilon,\xi_0+\epsilon]$. Denote by $\chi^+,\chi^-$
  the Euler characteristic of $P_{\xi_0+\epsilon},
  P_{\xi_0-\epsilon}$.  
  We have  
  \begin{enumerate}
  \item $\chi^+=\chi^-$, if $\xi_0$ is of type $E^{+}$ or $V^+$,
  \item $\chi^+=\chi^- - 1$, if $\xi_0$ is of type $E^{-}$ 
  \item $\chi^+=\chi^- - 2$, if $\xi_0$ is of type $V^-$, and 
  \item $\chi^+=\chi^- + 1$, if $\xi_0$ is of type $F^0$,
  \end{enumerate}
  Since $\chi(P_0) = \chi(\C^2) + \chi(S) - \chi(B^-\cap \C^2)$ 
  and $\chi(P_1)=  \chi(S) + \chi(B^+\cap \C^2)$, and since $H$
  has at most $\#(\C^0) = \chi(\C^0)$ critical points of type $V^-$,
  the lemma follows.  
\end{proof}

%***********************************************************************

\subsection{Relating triangulations by contractions and expansions} 
\label{sec:proof1}

This subsection is devoted to the proof of Theorem~\ref{thm:basic2}.
We start with a construction of a triangulation associated to embedded
fake surfaces, so that the triangulation changes by contractions
(resp.\ expansions) if the fake surface changes by deletion (resp.\ 
insertion) of $2$--cells.
The proof of Theorem~\ref{thm:basic2} is an application of this
construction to fake surfaces that are associated to the non-critical
parameters of a $\C^2$--Morse embedding. 

A fake surface $Q\subset S^3$ is \begriff{regular} if it is the
$2$--skeleton of a simple cellular decomposition $\C$ of $S^3$ such
that any $2$--cell of $\C$ is contained in the boundary of two
different $3$--cells of $\C$.
It follows easily that the closure of any open $3$--cell of $\C$ is a
compact ball, and the barycentric subdivision of $\C$ is a
triangulation of $S^3$ (i.e., it has no multiple edges).  Since $\C$
and its barycentric subdivision are determined by $Q\subset S^3$, we
denote this triangulation by $\T(Q)$. The next lemma provides
conditions under which the deletion of a $2$--stratum of $Q$ gives rise
to a sequence of contractions of $\T(Q)$.

\begin{lemma}\label{lem:TQ}
  Let $Q_1,Q_2\subset S^3$ be regular fake surfaces. Let $c$
  be a $2$--stratum of $Q_1$ whose closure contains $k$ intrinsic
  vertices. If $Q_2$ is obtained from $Q_1$ by deletion of $c$, then
  $\T(Q_2)$ is obtained from $\T(Q_1)$ by a sequence of $4k + 2$
  contractions.
\end{lemma}
\begin{proof}
  By hypothesis on $Q_1$, it is the $2$--skeleton of a simple cellular
  decomposition $\C_1$ of $S^3$, and $c$ is a $2$--cell contained in the
  boundary of two different $3$--cells $X_1,X_2$ of $\C_1$. 
  We contract $\T(Q_1)$ along the edges that connect the 
  barycenter of $c$ with the barycenters of $X_1,X_2$. By hypothesis 
  on $Q_2$, the closure of any connected component of $S^3\setminus Q$
  is a ball, hence $\d X_1\cap \d X_2 = \overline c$. Thus the two
  contractions are allowed, i.e., do not introduce multiple edges. 

  Any edge $e\subset \d c$ is adjacent to exactly two $2$--cells
  $c_1,c_2$ of $\C_1$ that are different from $c$. 
  By hypothesis on $Q_2$, we have $\d c_1\cap \d c_2=\overline e$. 
  Thus we can further contract along the 
  edges of $\T(Q_1)$ that connect the barycenter of $e$ with the
  barycenters of $c_1,c_2$, whithout to introduce multiple edges.

  Any vertex $v\in \d c$ is endpoint of exactly two edges $e_1,e_2$ of 
  $\C_1$ that are not contained in $\d c$. Since $\d e_1\cap \d e_2 = 
  v$, we can further contract along the edges of
  $\T(Q_1)$ that connect $v$ with the barycenters of $e_1,e_2$. 
  These $2+ 2k + 2k$ contractions yield $\T(Q_2)$. 
\end{proof}

Our plan is to associate fake surfaces to the non-critical parameters of 
a $\C^2$--Morse embedding, and apply to them the preceding techniques. 
For this aim, we need a $\C^2$--Morse embedding in a particularly nice
position, provided by the following lemma.

\begin{lemma}\label{lem:makeMorse}
  Let $\T$ be a triangulation of $S^3$, and let $\C$ be its dual
  cellular decomposition.
  There are two vertices $v_0,v_1\in \C^0$ and a $\C^2$--Morse
  embedding $H\co S^2\times I \to S^3$ with 
  $c(H,\C^1)\le  p(\T) + 2 n + 2$ and 
  $H(S^2\times I)= S^3\setminus U(\{v_0,v_1\})$, 
\end{lemma}
\begin{proof}
  By definition, there is a $\C^1$--Morse embedding $H'\co S^2\times
  I\to S^3$ with $\C^2\subset H'(S^2\times I)$ and
  $c(H', \C^1) = p(\T)$. Let $X_0,X_1$ be the open $3$--cells of $\C$
  that contain $H'_0,H'_1$. Pick two different  vertices $v_0,v_1$ in $\d
  X_0,\d X_1$. We change $H'$ into a $\C^1$--Morse embedding 
  $H''\co S^2\times I\to S^3$ by pushing a finger from $H'_0$ towards
  $v_0$ and and from $H'_1$ towards $v_1$, so that
  $\C^2\subset  H''(S^2\times I)$ and $v_0$ 
  (resp. $v_1$) are the critical points of $H''$ that correspond to
  the smallest (resp. biggest) critical parameter of $H''$ with respect
  to $\C^2$, denoted by
  $\xi_0$ (resp. $\xi_1$). 
  This is possible by introducing at  most one critical point in each
  of the eight edge germs at $v_0,v_1$. Thus $c(H'',\C^1) \le p(\T) +
  8$. 

  By a small isotopy, we factorize the critical points of $H''$ in
  $\C^2$, except $v_0,v_1$, by critical points of type $F^\pm$, $F^0$, 
  $E^\pm$ and $V^\pm$. 
  Any critical point of $H''$ in the interior of an edge of  $\C$
  factorizes by one critical point of type $E^\pm$ and some other
  critical points in $\C^2\setminus \C^1$. 
  Let $v_2\in \C^0\setminus\{v_0,v_1\}$ be a vertex, corresponding to
  a critical parameter $\xi_2$ of $H''$. If a component of
  $U(v_2)\setminus U(H''_{\xi_2})$ intersects $\C^1$ in exactly one
  (resp. two) arcs, then 
  the critical point $v_2$ factorizes by one (resp. two) critical
  points of type $E^\pm$, one critical point of type $V^\mp$ (with
  consistent signs), and some
  critical points outside $\C^1$. Thus we obtain a $\C^2$--Morse
  embedding $H'''\co S^2\times I\to S^3$ with $c(H''',\C^1)\le p(\T) +2 n + 
  4$. Now, a $\C^2$--Morse embedding $H\co S^2\times I\to S^3$ with the
  claimed properties is given by the restriction of $H'''$ to
  $S^2\times [\xi_0+\epsilon, \xi_1 - \epsilon]$, for small $\epsilon
  >0$.   
\end{proof}

To prove Theorem~\ref{thm:basic2}, it suffices to show that any
triangulation $\T$ of $S^3$ can be transformed into the barycentric
subdivision of the boundary complex of a $4$--simplex by $\le 325 p(\T)
+ 254$ expansions and contractions.
Let $\C$, $v_0$, $v_1$ and $H\co S^2\times I\to S^3$ be as in the
preceding lemma, with $H_0=\d U(v_0)$ and $H_1=\d U(v_1)$. For $\xi\in
I$, 
let $B^+(\xi)$ (resp.\ $B^-(\xi)$) be the closure of the component of
$S^3\setminus H_\xi$ that contains $v_1$ (resp.\ $v_0$).

Since $H_0\setminus\C^2$ and $H_1\setminus \C^2$ are disjoint unions
of discs, $H$ satisfies the general hypothesis of
Subsection~\ref{sec:minmax}. Hence by Proposition~\ref{lem:minweg}, we
can assume that $H$ has no critical points of type $F^\pm$.
By Proposition~\ref{lem:saddlebound}, $H$  has at most $c(H,\C^1)+2 n
+ 2$ critical points of type $F^0$ (hint: we have $\chi(B^-(0)\cap
\C^2)=\chi(U(v_{0})\cap \C^2)=1$, $\chi(B^+(1)\cap \C^2)=1$, and
$-\chi(\C^2) = \chi(S^3\setminus \C^2)\le n$).
With Lemma~\ref{lem:makeMorse}, we have $c(H,\C^2) \le 2 p(\T) + 6 n 
+ 6 \le 8 p(\T) + 6$.

We define $P_\xi = \left(\C^2\cap B^+(\xi)\right)\cup H_\xi$ for $\xi
\in I$. By the next lemma, the triangulations $\T_\xi=\T(P_\xi)$ of
$S^3$ are defined for any non-critical parameter $\xi\in I$ 
of $H$ with respect to $\C^2$.

\begin{lemma}\label{lem:Tdefined}
  For any non-critical parameter $\xi\in I$, $P_\xi$ is a regular fake surface.
\end{lemma}
\begin{proof}
  Let $c$ be a $2$--cell of $\C$ and assume that some component $\gamma$ 
  of $c\cap H_\xi$ is a circle. It bounds a disc $D\subset c$.
  Let a collar of $\gamma$ in $D$ be contained in $B^+(\xi)$ (resp.\ 
  in $B^-(\xi)$). Since $H$ has no critical parameters of type $F^-$ 
  (resp. $F^+$), it follows by induction on the number of critical
  parameters of $H$ that $c\cap H_1$ (resp.\ 
  $c\cap H_0$) contains a circle, in contradiction to the hypothesis
  on $H$.
  Thus $H_\xi$ intersects the $2$--cells of $\C$ in arcs.

  Similarly one shows that $H_\xi\setminus \C^2$ is 
  a disjoint union of open discs, and $S^3\setminus P_\xi$ is a 
  disjoint union of open $3$--balls. Thus the open $2$--strata of $P_\xi$
  are discs. 
  Since both $H_0$ and $H_1$ intersect $\C^1$ and $\C^1$ is connected,
  $P_\xi$ has an intrinsic vertex in $H_\xi\cap \C^1$.
  Since $\xi$ is not a critical parameter of $H$, 
  $P_\xi$ is simple.
  In conclusion, $P_\xi$ is the $2$--skeleton of a simple cellular decomposition
  $\C_\xi$ of $S^3$.

  Let $c$ be a $2$--cell of $\C_\xi$. If $c\subset \C^2$ then it separates
  two $3$--cells of $\C_\xi$, since $\C$ is dual to a triangulation.
  If $c\subset H_\xi$ then it separates the $3$--cell of $\C_\xi$ 
  corresponding to $B^-(\xi)$ from another $3$--cell of $\C_\xi$. 
\end{proof}

We show how $P_\xi$ and $\T_\xi$ change when $\xi$ passes a critical
parameter $\xi_0$ of $H$ with respect to $\C^2$. Let
$p_0\in \C^2$ be the critical point corresponding to $\xi_0$. 
We choose $\epsilon > 0$ so  that $\xi_0$ is the only  
critical parameter in $[\xi_0-\epsilon,\xi_0+\epsilon]$. Choose
local coordinates $(x,y,z)$ around $p_0$ as in
Definition~\ref{def:types}. Let $r>0$ be small enough such that
$B=\{x^2+y^2+z^2 \le r^2\}$ is a closed regular neighborhood of $p_0$. 

By isotopy of $H_{\xi_0\pm\epsilon}$ mod $\C^2$, we can assume that $B
\cap H_{\xi_0-\epsilon}=D$ and $B\cap H_{\xi_0+\epsilon}=D'$ are
discs, $\d B= D\cup D'$, and $H_{\xi_0+\epsilon}=
(H_{\xi_0-\epsilon}\setminus D)\cup D'$. We define $P'=
(P_{\xi_0-\epsilon}\setminus B)\cup \d B$. One easily verifies that
$P'$ is a regular fake surface, hence $\T(P')$ is defined.
By deletion of its $2$--stratum $D$, one obtains $P_{\xi_0+\epsilon}$. In
$\d D$ are at most $4$ intrinsic vertices (namely when $p_0$ is of type
$F^0$). Thus by Lemma~\ref{lem:TQ}, $\T_{\xi_0+\epsilon}$ is obtained
from $\T(P')$ by $\le 18$ contractions.
We consider how $P'$ and$P_{\xi_0-\epsilon}$  are related by
deletions of $2$--strata. 
\begin{enumerate}
\item If $p_0$ is of type $F^0$, then one obtains $P_{\xi_0-\epsilon}$ 
  from $P'$ up to isotopy by deletion of the two $2$--strata corresponding
  to $D'\cap \{z\le x^2-y^2\}$. They both have $2$ intrinsic vertices
  in its boundary.
\item If $p_0$ is of type $E^{+}$, then one obtains
  $P_{\xi_0-\epsilon}$ from $P'$ up to isotopy by
  deletion of the $2$--stratum corresponding to $D'\cap \{z\le x^2 + y, 
  y\ge 0\}$. It has $4$ intrinsic vertices in its boundary.
\item If $p_0$ is of type $E^{-}$, then $P_{\xi_0-\epsilon}$ is
  isotopic to $P'$. 
\item If $p_0$ is of type $V^+$, then one obtains $P_{\xi_0-\epsilon}$ 
  from $P'$ up to isotopy by deletion of the $2$--stratum corresponding to
  $D'\cap \{z\le |x| + y, y\ge 0\}$. It has $5$ intrinsic vertices in its
  boundary. 
\item If $p_0$ is of type $V^-$, then one obtains $P_{\xi_0-\epsilon}$ 
  from $P'$ up to isotopy by insertion of the $2$--stratum 
  corresponding to $B\cap \{z\le -|y|, x=0\}$. It has $3$ intrinsic vertices in 
  its boundary.
\end{enumerate}
Thus by Lemma~\ref{lem:TQ}, $\T(P')$ is obtained from
$\T_{\xi_0-\epsilon}$ by $\le 22$ successive contractions or
expansions.

With our bound for $c(H,\C^2)$, it follows
that $\T_0$ and $\T_1$ are related by $$\le (18 + 22)\cdot \left(8
  p(\T) + 6\right) = 320 p(\T) + 240$$  expansions and
contractions. 
Let $n$ be the number of tetrahedra of $\T$. Since $\T$ has $2 n$
$2$--simplices and at most $2 n$ edges, one obtains its barycentric
subdivision $\T(\C^2)$ by $\le 5 n \le 5 p(\T)$ 
expansions.
Since $Q_0$ is isotopic to the result of adding one triangular $2$--stratum
to $\C^2$, $\T(\C^2)$ can be transformed into $\T_0$ by $14$
expansions. 
We count together and find that $\T$ can be transformed into $\T_1$,
the barycentric subdivision of the boundary complex of a $4$--simplex, by 
$\le 5 p(\T) + 14 + 320 p(\T) + 240 = 325 p(\T) + 254$ 
contractions and  expansions. This yields 
Theorem~\ref{thm:basic2}.

%************************************************************

\subsection{How to make a triangulation edge contractible} 
\label{sec:proof2}

In this subsection, we prove Theorem~\ref{thm:basic3}.
We use the same notations as in the previous subsection.
The idea for the proof is to insert
$2$--strata as above in the transition from $P_{\xi_0-\epsilon}$ to $P'$,
and to postpone the deletion of $2$--strata until all
insertions are done. 
This means to transform $\T$ by expansions into
a triangulation that is then turned into the boundary complex
of a $4$--simplex by a sequence of contractions.

Let $\xi_1<\xi_2< \dots < \xi_Z$ be the critical parameters of $H$ with
respect to $\C^2$, let $p_1,\dots,p_Z\in \C^2$ be the corresponding
critical points, and let $\xi_{Z+1}=1$.
Let $\epsilon > 0$ be small enough such that $\xi_i$ is the only critical parameter
of $H$ in $[\xi_i-\epsilon,\xi_i+\epsilon]$, for all
$i=1,\dots,Z$. 
For any $i$, choose local coordinates $(x_i,y_i,z_i)$ around $p_i$ as 
in Definition~\ref{def:types}. 
Let $r>0$ be small enough such that $B_i=\{x_i^2+y_i^2+z_i^2 \le r^2\}$ is a closed
regular neighborhood of $p_i$. 
We arrange $H_{\xi_i\pm\epsilon}$ by isotopy mod $\C^2$ so that 
$\d B_i \cap H_{\xi_i-\epsilon}=D_i$ and  $\d B_i\cap 
H_{\xi_i+\epsilon}=D'_i$ are discs, $\d B_i= D_i\cup D'_i$, and
$H_{\xi_i+\epsilon}= (H_{\xi_i-\epsilon}\setminus D_i)\cup D'_i$.

We now define iteratively a sequence $Q_1,\ldots,Q_{Z+1}\subset S^3$
of regular fake surfaces, together with graphs $\Gamma_i\subset
Q_i^1$, so that $P_{\xi_i-\epsilon}\subset Q_i$ and $\Gamma_i
= \d (Q_i\setminus P_{\xi_i-\epsilon})$.
We define $Q_1 = (\C^2\cap B^+(\xi_1-\epsilon))\cup H_{\xi_1-\epsilon}= P_{\xi_1-\epsilon}$ 
and $\Gamma_1=\emptyset$. 
Let $i\in\{1,\dots, Z\}$. We arrange $\Gamma_i\subset Q_i^1$ by an
isotopy of $Q_i$ mod $\C^2$ so that 
$\Gamma_i$ intersects $\d D_i$ transversely and $\#(\Gamma_i\cap \d
D_i)$ is as small as possible.  
If $\xi_i$ is not of type $V^-$, then we define $Q'_i = (Q_i \setminus
B_i)\cup \d B_i$ and $\Gamma'_i = (\Gamma_i\setminus D_i)\cup \d D_i$.
If $\xi_i$ is of type $V^-$, then we define $Q'_i = (Q_i \setminus B_i) \cup
\d B_i \cup (\{x_i=0\}\cap B_i)$ and $\Gamma'_i = (\Gamma_i\setminus
D_i)\cup \d D_i \cup (\{x_i=0\}\cap D'_i)$; here $Q'$ is isotopic to $Q_i$.
We define $Q_{i+1}$ and $\Gamma_{i+1}$ as the result of $Q'_{i}$ and
$\Gamma'_i$ under an isotopy mod $\C^2$ that relates
$H_{\xi_i+\epsilon}$ with $H_{\xi_{i+1}-\epsilon}$.  

It follows as in the proof of Lemma~\ref{lem:Tdefined} that $Q_1,\dots
,Q_{Z+1}$ are regular fake surfaces. Let $\T_i=\T(Q_i)$, for
$i=1,\dots,Z+1$. 
In the following two lemmas, we show that for $i=1,\dots, Z$ one
obtains $Q_i$ from $Q_{i+1}$ by deletion of $2$--strata, with an
estimate for the number of vertices in the boundary of these $2$--strata.
\begin{lemma}\label{lem:countGamma}
  For $i=1,\dots,Z$, we have $\#(\Gamma_i\cap \C^2) \le 4(i-1)$.
\end{lemma}
\begin{proof}
  Let $j\ge 1$.
  One observes that $\Gamma'_{j}\cap \C^2$ comprises at most four
  points more than $\Gamma_j\cap \C^2$, namely the points of $\d
  D_i\cap \C^2$. The claim follows by induction, with
  $\Gamma_1=\emptyset$.    
\end{proof}
\begin{lemma}\label{lem:countcontr}
  For $i=1,\dots,Z$, one obtains $Q_i$ from $Q_{i+1}$ up to isotopy
  by deletion of one $2$--stratum 
  with at most $5$ vertices in its boundary, or by deletion of two
  $2$--strata of $Q_{i+1}$ with at most 
  $2i$ vertices in its boundary.
\end{lemma}
\begin{proof}
  Let $p_i$ be not of type $F^0$. Then any simple arc in 
  $D_i\setminus \C^2$ with boundary in $\d D_i$ 
  is parallel in $D_i\setminus \C^2$ to a sub-arc of $\d D_i\setminus
  \C^2$. Thus $U(p_i)\cap \Gamma_i=\emptyset$ by the minimality of
  $\#(\Gamma_i\cap \d D_i)$. If
  $p_i$ is of type $V^-$, then $Q_i$ is isotopic to $Q_{i+1}$. If
  $p_i$ is of type $E^\pm$ or $V^+$, then we 
  apply the analysis in the proof of Theorem~\ref{thm:basic2} and see
  that one obtains $Q_i$ from $Q_i'\simeq Q_{i+1}$ by deletion of one
  $2$--stratum with at most $5$ vertices in the boundary. 

  If $p_i$ is of type $F^0$, then one obtains $Q_i$ from $Q'_i\simeq
  Q_{i+1}$ by deletion of the two connected components of $D'_i\cap
  \{z_i\le x_i^2-y_i^2\}$. 
  We estimate the number of intrinsic vertices in the boundary of these
  $2$--strata. 
  Let $c$ be the $2$--stratum of $Q_i$ that contains
  $D_i\cap  \{z_i\le x_i^2-y_i^2\}$. 
  By minimality of $\Gamma_i\cap \d D_i$, there are at most $\frac 12
  \#(\Gamma_i\cap \d c)$ arcs in $\Gamma_i\cap D_i$, each connecting the
  two components of $c\cap \d D_i$.  
  Thus by Lemma~\ref{lem:countGamma} both components  of $D'_i\cap
  \{z_i\le x_i^2-y_i^2\}$ have at
  most $2(i-1) + 2$ vertices in their boundary. 
\end{proof}

By Lemmas~\ref{lem:TQ} and~\ref{lem:countcontr}, one obtains
$\T_{i+1}$ from $\T_{i}$ by at most $\max\{22, 16i + 4\}$
expansions. As in the previous section, $\T_1$ is obtained from $\T$ 
by $\le 5p(\T) + 14$  expansions, and $Z\le 8p(\T) + 6$. Thus one 
obtains $\T_{Z+1}$ from $\T$ by  
\begin{eqnarray*} 
  &\le&  5p(\T) + 14 + 22 + \sum_{i=2}^Z (16 i + 4) \\
  &\le&     512 (p(\T))^2 + 869 p(\T) + 376
\end{eqnarray*}
expansions.

One obtains $(\C^2\cap B^+(1-\epsilon))\cup
H_{1-\epsilon}=P_{1-\epsilon}$ from $Q_{Z+1}$ by successive 
deletion of the $2$--strata $D_1, D_2,\dots, D_Z$. One checks that these
deletions satisfy the hypo\-the\-sis of Lemma~\ref{lem:TQ}. Thus
$\T(P_{1-\epsilon})$ is the result of $\T_{Z+1}$ under successive
contractions.
Since $H_{\xi_Z+\epsilon}$ is isotopic mod $\C^2$ to $\d U(v_1)$,
$\T(P_{\xi_Z+\epsilon})$ is the barycentric subdivision of the
boundary complex of a $4$--simplex. Thus $\T_{Z+1}$ is edge contractible,
which finally proves Theorem~\ref{thm:basic3}.

%************************************************************************
%************************************************************************

\begin{center}\textbf{Acknowledgments}\end{center}

This paper is part of my PhD thesis, defended in June 2001 at the
Institut de Recherche Math\'ematique Avanc\'ee in Strasbourg/France.
I am grateful to my supervisor Vladimir G. Turaev for many fruitful
conversations. I was supported by a grant of the Studienstiftung des
deutschen Volkes.

%************************************************************************
%************************************************************************

\end{document}